\documentclass[10pt,twoside]{article}
\usepackage{bcc24book,url,booktabs,graphicx}
\usepackage{eucal,amsfonts,amsxtra}
\usepackage[nobysame]{amsrefs}

    \newcommand{\LO}{{\mathcal{LO}}}
 
    \newcommand{\bA}{\mathfrak A}
    \newcommand{\bB}{\mathfrak B}
    \newcommand{\bF}{\mathfrak F}
    \newcommand{\bD}{\mathfrak D}
    \newcommand{\bE}{\mathfrak E}
        \newcommand{\bC}{\mathfrak C}
        \newcommand{\bX}{\mathfrak X}
       \newcommand{\bY}{\mathfrak Y}
       
  \newcommand{\bT}{\mathfrak{T}}
    \newcommand{\bP}{\mathfrak{P}}

\newcommand{\C}{{\mathcal C}}

\DeclareMathOperator{\Age}{Age}
\DeclareMathOperator{\Aut}{Aut}

\DeclareMathOperator{\Betw}{Betw}
\DeclareMathOperator{\Forb}{Forb}

\newcommand{\mN}{\mathbb N}
\newcommand{\mQ}{\mathbb Q}
\newcommand{\mV}{\mathbb V}

\newcommand{\qed}{$\square$}

\newcommand{\Fresse}{Fra\"{i}ss\'{e}}
    \newcommand{\Nesetril}{Ne\v{s}et\v{r}il}

\newcommand{\ignore}[1]{}

\bcctitle{Ramsey Classes: Examples and Constructions}
\bccshorttitle{Ramsey Classes}

\bccname{Manuel Bodirsky\footnote{The author has received funding from the European Research Council under the European Community's Seventh Framework Programme (FP7/2007-2013 Grant Agreement no. 257039).}}
\bccshortname{M.~Bodirsky}
   \bccaddressa{Institut f\"ur Algebra, TU Dresden, \\ 01069 Dresden, Germany}
    \bccemaila{manuel.bodirsky@tu-dresden.de}

\begin{document}

\makebcctitle

\begin{abstract}
This article is concerned
with classes of relational structures 
that are closed under taking
substructures and isomorphism, that have the joint embedding property, and that furthermore have the \emph{Ramsey property}, a strong 
combinatorial property which resembles the 
statement of Ramsey's classic theorem. 
Such classes of structures have been called
\emph{Ramsey classes}. %; they
%Let $\mathcal C$ be a class of structures
%that is closed under taking substructures and that has the joint embedding property. 
% Nicht nur relational wegen Konstanten.
%Then $\mathcal C$ is said to be a \emph{Ramsey class}
%if for all $A,B \in \C$ there exists $C \in \cal C$ such
%that for every 2-colouring of the copies of $A$ in $C$ there exists a copy $B'$ of
%$B$ in $C$ such that all copies of $A$ in $B'$ 
%have the same colour. 
%have been introduced by \Nesetril, 
\Nesetril\ and R\"odl showed that 
they 
have the \emph{amalgamation property}, and 
therefore each such class has a homogeneous \Fresse\ limit. Ramsey classes have recently attracted attention due to a surprising link with
the notion of extreme amenability from topological dynamics.
Other applications of Ramsey classes include reduct
classification of homogeneous structures. 

We give a survey of the various fundamental
Ramsey classes and their (often tricky) combinatorial proofs,
and about various %(less complicated) 
methods to derive new Ramsey
classes from known Ramsey classes. Finally, we state open problems related to a potential classification of Ramsey classes. 
\end{abstract}

\tableofcontents
\thispagestyle{empty}
\setcounter{page}{1}  

\section{Introduction}
Let $\mathcal C$ be a class of finite relational structures.
Then $\mathcal C$ has the \emph{Ramsey property} if it satisfies a property that resembles 
the statement of Ramsey's theorem: 
%~\cite{Ram:On-a-problem}: 
for all $\bA,\bB \in \C$ there exists $\bC \in \cal C$ such that for every colouring of the embeddings of $\bA$ into $\bC$ with finitely many colours 
there exists a `monochromatic copy' of $\bB$ in $\bC$, that is, an embedding $e$ of $\bB$ into $\bC$ such that all embeddings of $\bA$ into the image of $e$ have the same colour. 
An example of a class of structures with the Ramsey property is the class of all finite
linearly ordered sets; this is Ramsey's theorem~\cite{Ram:On-a-problem}. Another example of a class with the Ramsey property is
the class of all ordered finite graphs, that is, 
structures $(V;E,\preceq)$ where $V$ is a finite set, $E$
the undirected edge relation, and $\preceq$ a linear order on $V$;
this result has been discovered by \Nesetril\ and R\"odl~\cite{NesetrilRoedl}, and, independently, Abramson and Harrington~\cite{AbramsonHarrington}. 

In this article we will be concerned exclusively with
classes $\mathcal C$ that are 
closed under taking substructures and isomorphism, and that have the \emph{joint embedding property}:
whenever $\bA,\bB \in \mathcal \mathcal C$,
then there exists a $\bC \in \mathcal C$ 
such that both $\bA$ and $\bB$ embed into $\bC$.
These are precisely the classes $\mathcal C$
for which there exists a countably infinite structure $\Gamma$ such that a structure belongs to
$\mathcal C$ if and only if it embeds into $\Gamma$. 
This statement also holds when the relational
signature of $\mathcal C$ is infinite, but here we additionally require that the class $\mathcal C$
has only countably many non-isomorphic members.
Following \Fresse's terminology, we say
that $\mathcal C$ is the \emph{age} of $\Gamma$. 

A class $\mathcal C$ will be called a \emph{Ramsey class}~\cite{RamseyClasses} if it has the Ramsey property, and is the age of a countable structure. 
%is closed under isomorphisms and substructures, and has the joint embedding property.
It is an open research problem, raised in~\cite{RamseyClasses}, whether Ramsey
classes can be \emph{classified} in some sense that needs to be specified. 

It has been shown by \Nesetril~\cite{RamseyClasses} that
Ramsey classes have the \emph{amalgamation property}, a central property in model theory. 
A class of structures $\mathcal C$ has the amalgamation property 
if for all $\bA,\bB_1,\bB_2 \in \mathcal C$
with embeddings $e_i$ 
of $\bA$ into $\bB_i$, for $i \in \{1,2\}$, there exist $\bC \in \mathcal C$ and embeddings $f_i$ of $\bB_i$ into $\bC$ such
that $f_1(e_1(a)) = f_2(e_2(a))$ for all elements $a$ of $\bA$.
A class of finite relational structures $\mathcal C$ is an \emph{amalgamation class} if it is closed under induced substructures, isomorphism, has countably many non-isomorphic members, and
the amalgamation property.
By \Fresse's theorem (which will be recalled in Section~\ref{sect:amalgamation})
%(see e.g.~\cite{Hodges}),
for every amalgamation 
class $\mathcal C$ 
there exists a countably infinite structure $\Gamma$ 
of age $\mathcal C$
%, such that 
%\begin{itemize}
%\item the \emph{age} of $\Gamma$, that is, the class of finite structures that embeds into $\Gamma$, equals $\mathcal C$; 
%\item 
which is \emph{homogeneous}, that is, any isomorphism between finite substructures of $\Gamma$ can be extended to an automorphism of $\Gamma$. 
The structure $\Gamma$ is in fact unique up to isomorphism, and called the \emph{\Fresse\ limit of $\mathcal C$}. In our example above 
where $\mathcal C$ is the class of all finite linearly ordered sets $(V;<)$, the \Fresse\ limit is isomorphic to $({\mathbb Q};<)$, that is, the
linear order of the rationals. 
%\end{itemize}

The age of a homogeneous structure with a finite relational signature is in general not Ramsey. However, quite surprisingly, homogeneous structures with finite relational signature typically have a homogeneous \emph{expansion}
by finitely many relations such that the
age of the resulting structure is Ramsey. 
The question whether we can replace in the previous sentence the
word `typically' by the word `always' appeared
in discussions of the author with Michael Pinsker and Todor Tsankov in 2010, and
has been asked, first implicitly in a conference publication~\cite{BPT-decidability-of-definability}, then explicitly in the journal version. 
The question 
motivates much
of the material present in this article, so we prominently state it here as follows. 

\begin{conj}[Ramsey expansion conjecture]
\label{conj:ramsey-expansion}
%\begin{itemize}
%\item 
Let $\Gamma$ be a homogeneous structure with finite relational signature. Then $\Gamma$ has a homogeneous  expansion  by finitely many relations whose age has the Ramsey property.  
%\item Let $\Gamma$ be $\omega$-categorical.
%Then $\Gamma$ has an $\omega$-categorical expansion $\Delta$ such that $\Delta$ is $\omega$-categorical and Ramsey. 
%\end{itemize}
\end{conj}
This conjecture has explicitly been confirmed
for all countable 
homogeneous directed graphs in~\cite{LafJasTheWoo} (those graphs have been
classified by Cherlin~\cite{Cherlin}),
and other homogeneous structures of interest~\cite{Hubicka-Nesetril-Bowtie}. 
%
%The conjecture is clearly related to the question whether Ramsey classes can be classified: 
%a classification
%of Ramsey classes probably helps to verify or falsify the conjecture; conversely, a proof of the conjecture  might 
%would be detailed enough
% should probably 
%allow to verify or falsify the conjecture 
%
%hich is an interesting example that has algebraicity~\cite{
%an $\omega$-categorical structure with algebraicity
%The conjecture is also true for all
%structures that are additionally $\omega$-stable  (Proposition~\ref{prop:}). 
The Ramsey expansion conjecture has 
several variants that are formally unrelated, but related in spirit; we will come back to this in the final section of the article. There we also discuss that the conjecture can be
translated into questions in topological dynamics which
are of independent interest. 

This text has its
focus on the combinatorial aspects of the theory, rather
than the links with topological dynamics. 
%Demonstrated in this text: we avoid topological dynamics entirely in this article (and in this way spare the reader of many definitions). 
What we do find convenient, though,
is the usage of concepts from model theory to present the results: instead of manipulating amalgamation classes $\mathcal C$
% of structures $\mathcal C$ that are closed under substructures, isomorphism, and that have the joint embedding and the amalgamation property, 
it is often more convenient to directly manipulate the homogeneous structures of age $\mathcal C$.

%also appears in~\cite{MelTheTsa,LafJasTheWoo}, and 

%\subsection*{Outline of the article}
\paragraph{Outline of the article.}
In Section~\ref{sect:ramsey-classes}
we give a self-contained introduction
to the basics of Ramsey classes, including
the proofs of some well-known and easy observations about them. In Section~\ref{sect:new-old} we show
how to derive new Ramsey classes from known ones; this section contains various facts or proofs that have not explicitly appeared in the literature yet.
\begin{itemize}
\item In %Section~\ref{sect:products} and~
Section~\ref{sect:interpret} we have basic results about the Ramsey properties of interpreted structures that have not been formulated previously in this form, but 
that are not difficult to show via variations of the so-called \emph{product Ramsey theorem}. 
\item In Section~\ref{sect:constants} we present a new non-topological proof,  due to Miodrag Sokic, of a known fact from~\cite{BPT-decidability-of-definability} about expanding Ramsey classes with constants.
\item In Section~\ref{sect:mc} and~\ref{sect:mc-core} we present generalisations of results from~\cite{Bod-New-Ramsey-classes}
about the Ramsey properties of model-companions and model-complete cores
of $\omega$-categorical structures. 
\end{itemize}

Some fundamental Ramsey classes
cannot be constructed by the general construction principles from Section~\ref{sect:new-old}. The most powerful tool
that we have to prove Ramsey theorems
from scratch is the \emph{partite method}, developed in the 70s and 80s, most notably by \Nesetril\ and R\"odl, 
which we present in Section~\ref{sect:partite}.
With this method we will show that the following
classes are Ramsey: the class of all ordered graphs, the class of all ordered triangle-free graphs, or more generally the class of all ordered structures given by a set of homomorphically forbidden irreducible substructures. 
%Finally, we generalise this
%even further, and prove that classes of ordered
%structures that are given by homomorphically forbidding finitely many finite structures 
%have a finite homogeneous Ramsey expansion,
%and therefore provide further examples for Conjecture~\ref{conj:ramsey-expansion}. 

%At the workshop on homogeneous structures at %the university of Leeds in July 2011, 
%\Nesetril\ has asked whether for all Ramsey classes that Ramsey property can be shown with some variation of the partite method. 

There are also Ramsey classes with finite relational signature where it is not clear how to show the 
Ramsey property with the partite method, to the best of my knowledge.
We will see such an example, based on Ramsey theorems for tree-like structures, in 
Section~\ref{sect:inductive}.

When we want to make progress on Conjecture~\ref{conj:ramsey-expansion}, we need a better understanding of the type of expansion
needed to turn a homogeneous structure in a finite language into a Ramsey structure. 
Very often, this can be done by adding a
linear ordering to the signature (a partial explanation for this is given in Section~\ref{sect:omega-cat}). 
But not any linear ordering might do the job;
a crucial property for finding the right 
ordering is the so-called \emph{ordering property},
which is a classical notion in structural Ramsey theory. 
We will present in Section~\ref{sect:ordering} 
a powerful condition that implies
that a Ramsey class has the ordering property
with respect to some given ordering.
%and a recent result of Zucker that shows
%that proving Conjecture~\ref{conj:ramsey-expansion}
%is equivalent to proving finiteness of Ramsey degrees.

Finally, in Section~\ref{sect:open}, we discuss the mentioned link between Ramsey theory and topological dynamics, then present an application 
of Ramsey theory for classifying reducts of homogeneous structures, and conclude 
with some open problems 
related to Conjecture~\ref{conj:ramsey-expansion}. 

\section{Ramsey classes: definition, examples, background}
The definition of Ramsey classes is inspired 
by the statement of the classic
theorem of Ramsey, which we therefore 
recall in the next subsection, before defining
the Ramsey property in Section~\ref{sect:ramsey-classes}
and Ramsey classes in Section~\ref{sect:jep}. 

There are two important necessary conditions for
a class to be Ramsey: \emph{rigidity} (Section~\ref{sect:rigidity})
and \emph{amalgamation} (Section~\ref{sect:amalgamation}). 
We will see examples that show that these two
conditions are not sufficient (Section~\ref{sect:counter}). 
The Ramsey property of a Ramsey class $\mathcal C$
can be seen as a property of
the automorphism group of the \Fresse\ limit of $\mathcal C$; this perspective is discussed in Sections~\ref{sect:automorphisms} and~\ref{sect:omega-cat}.

% JEP
% Rigidity
% Amalgamation
% Counter-examples
% Automorphism groups
% omega-categoricity

\subsection{Ramsey's theorem}
The set of positive integers is denoted by $\mN$, and the set $\{1,\dots,n\}$ is denoted by $[n]$.
%Subsets of a set of cardinality $m$ will be called $m$-subsets in the following. 
For $M,S \subseteq \mN$ we write
${M \choose S}$ for the set of all order-preserving maps from $S$ into $M$. 
When $f$ is a map, and $\mathcal S$ is a set of maps 
whose range equals the domain of $f$, then
$f \circ \mathcal S$ denotes the set
$\{f \circ e \mid e \in \mathcal S\}$. 
%The expression $f \circ {M \choose S}$
%denotes $\{f \circ e \mid e \in {M \choose S}\}$. 
%We also refer to mappings $f \colon {S \choose m} \rightarrow [r]$ as a \emph{colouring}
%of $S$ (with the \emph{colours} $[r]$). 
A proof of Ramsey's theorem can be found
in almost any textbook on combinatorics. 

\begin{theorem}[Ramsey's theorem~\cite{Ram:On-a-problem}]\label{thm:ramsey}
For all $r,m,k \in \mN$ there is a positive integer $g$ such that for every $\chi \colon {[g] \choose [k]} \rightarrow [r]$ there exists an $f \in {[g] \choose [m]}$ such that $|\chi(f \circ {[m] \choose [k]})| \leq 1$.  
\end{theorem}

\subsection{The Ramsey property}
\label{sect:ramsey-classes}
In this section we 
define the Ramsey property for
classes of structures.
All structures in this article have an at most countable 
domain, and have an at most countable signature. Typically, the signature will be relational
and even finite; but many results generalise to
signatures that are infinite and also contain function symbols. 
In Section~\ref{sect:constants} 
%and Section~\ref{sect:canonical}, 
it will be useful to consider signatures that also contain constant symbols (i.e., function symbols of arity zero).

Let $\tau$ be a relational 
signature, let $\bB$ be a $\tau$-structure.
For $R \in \tau$, we write $R^\bB$ for the corresponding relation of $\bB$. Typically, 
the domain of $\bA,\bB,\bC$ will be denoted by $A,B,C$, respectively.  
Let $A$ be a subset of the domain $B$ of $\bB$. 
Then the \emph{substructure} of $\bB$ \emph{induced
by $A$} is the $\tau$-structure $\bA$ with domain
$A$ such that for every relation symbol $R \in \tau$ of arity $k$ 
we have $R^\bA = R^\bB \cap A^k$.

If $\tau$ is not a purely relational signature, but
also contains constant
symbols, then every substructure $\bA$ of 
$\bB$ 
must contain for every constant symbol $c$ in $\tau$ 
the element $c^\bB$, and 
$c^\bA = c^\bB$.  
%We say that $A$ is a \emph{substructure}
%of $B$ if the domain $D$ of $A$
%is a subset of the domain of $B$, 
%if for every relation symbol $R \in \tau$ of arity $k$ we have $R^A = R^B \cap D^k$.
%, and for every function symbol $f \in \tau$ of arity $k$ and $d_1,\dots,d_k \in D$ we have that 
%$f^A(d_1,\dots,d_k) = f^B(d_1,\dots,d_k)$. In particular, when $c \in \tau$ is a constant  symbol (a relation symbol of arity zero), then $c^B$ is an element of all substructures of $B$.
%For any subset $S$ of the domain of $B$, the \emph{substructure induced by $S$ in $B$} is the substructure of $B$ whose domain contains $A$ and is minimal with this property. 
An \emph{embedding} of $\bB$ into $\bA$ is a
mapping $f$ from $B$ to $A$ which is an isomorphism between $\bB$ and the substructure induced by the image of $f$ in $\bB$.
This substructure will also be called a \emph{copy} of $\bA$ in $\bB$. 
%, and is denoted by $f[A]$. 
We write ${\bB \choose \bA}$ for the set of all 
embeddings of $\bA$ into $\bB$.

\begin{Def}[The partition arrow]
When $\bA,\bB,\bC$ are $\tau$-structures, and $r \in \mathbb N$, then we write
$\bC \to (\bB)^{\bA}_r$
if for all $\chi \colon {\bC \choose \bA} \to [r]$ there exists an 
$f \in {\bC \choose \bB}$ such that 
$|\chi(f \circ {\bB \choose \bA})| \leq 1$. 
\end{Def}

We would like to mention that in some 
papers, the partition arrow is defined 
for the situation where ${\bB \choose \bA}$ does
not denote the set of embeddings of $\bA$ into $\bB$, but the set of copies of $\bA$ in $\bB$. 
These two definitions are closely related; the article~\cite{MuellerPongracz} is specifically 
about this difference. 
Also~\cite{Topo-Dynamics} and~\cite{Zucker14} treat the relationship between the two definitions.

In analogy to the statement of Ramsey's theorem, we can now define the Ramsey property
for a class of relational structures.

\begin{Def}[The Ramsey property]
\label{def:Ramseyclass}
A class $\mathcal C$ of finite structures has the 
\emph{Ramsey property} if for all $\bA,\bB \in \mathcal C$  and  $k \in \mN$ there exists a $\bC \in \mathcal C$
such that $\bC \to (\bB)^{\bA}_k$. 
\end{Def}

\begin{example}\label{expl:lo}
The class of all finite linear orders, denoted by $\LO$, has the Ramsey property. This is a reformulation of Theorem~\ref{thm:ramsey}. \qed
\end{example}

The following well-known fact shows that we can
always work with 2-colourings instead of general
colourings when we want to prove that a certain class has the Ramsey property. 

\begin{lemma}\label{lem:colours}
Let $\mathcal C$ be a class of structures,
and $\bA \in \mathcal C$. Then for every $\bB \in \mathcal C$ and $r \in \mN$ there exists a $\bC \in \mathcal C$ such that $\bC \to (\bB)^\bA_r$
if and only if for every $\bB \in \mathcal C$ there exists a $\bC \in \mathcal C$ such that $\bC \to (\bB)^\bA_2$. 
\end{lemma}
\begin{proof}
Suppose that for every $\bB \in \mathcal C$ there exists a $\bC \in \mathcal C$ such that $\bC \to (\bB)^\bA_2$. We inductively define a sequence $\bC_1,\dots,\bC_{r-1}$ of structures in $\mathcal C$ as follows. 
Let $\bC_1$ be such that $\bC_1 \to (\bB)^\bA_2$.
For $i \in \{2,\dots,r-1\}$, let $\bC_i$ be such that
$\bC_i \to (\bC_{i-1})^\bA_2$. We leave it to the reader to verify that $\bC_{r-1} \to (\bB)^\bA_r$. 
%Now, let $\chi \colon {C_k \choose A} \to [k]$ be arbitrary. For $i \in \{1,\dots,k\}$, 
%define $\xi_i \colon {C_k \choose A} \to [2]$ 
%as follows. For $e \in {C_k \choose A}$, 
%define $\xi_i(e) = 1$ if $\chi(e) = i$ and
%$\xi_1(e) = 2$ otherwise. 
%We now inductively define $e_{k-1}$, \dots, $e_1$, where $e_i$ is an embedding of
%$C_i$ into $C_{i+1}$. 
%Since 
%$C_k \to (C_{k-1})^A_2$ there is an $e \in {C_k \choose C_{k-1}}$ such that
%$\xi_k$ is monochromatic on ${e_{k-1}[C_{k-1}] \choose A}$. If $\xi_k$ is constant $1$ then $\chi$ is constant on ${e_{k-1}[C_{k-1}] \choose A}$;
%since there is an embedding of $C$ into $C_{k+1}$ we are done. Otherwise, 
%$\chi({e_{k-1}[C_{k-1}] \choose A}) \subseteq \{1,\dots,k-1\}$.  Set $e_{k-1} := e$. 
%Suppose we have already defined $e_i$,
%and that we want to define $e_{i-1}$. 
%Since $C_k \to (C_{i-1})^A_2$ there is an $e_{i-1} \in {e_i[C_i] \choose C_{i-1}}$ such that
%$\xi_i$ is monochromatic on ${e_{i-1}[C_{i-1}] \choose A}$. We again have that either
%$\xi_i$ is constant $1$ and we are done, or 
%$\chi({e_{i-1}[C_{i-1}] \choose A}) \subseteq \{1,\dots,i-1\})$. In particular, we have 
%$\chi({{e_{1}[C_1] \choose A}) = \{1\}$,
%and we have found a monochromatic copy
\end{proof}

\subsection{The joint embedding property and Ramsey classes}
\label{sect:jep}
We say that a class of structures $\mathcal C$
is \emph{closed under substructures}
if for every $\bB \in \mathcal C$, all substructures of $\bB$ are also in $\mathcal C$. 
The class $\mathcal C$
is \emph{closed under isomorphism}
if for every $\bB \in \mathcal C$, all structures that are isomorphic to $\bB$ are also in $\mathcal C$. 
In this article, we will focus on classes of finite
structures that are closed under induced substructures and isomorphism, and that have the joint embedding property. Recall from the introduction that $\mathcal C$ has the joint
embedding property
if for every $\bA,\bB
\in \mathcal \mathcal C$,
there exists a $\bC \in \mathcal C$ 
such that both $\bA$ and $\bB$ embed into $\bC$.
Such classes of structures naturally arise
as follows; see e.g.~\cite{Hodges}. 

\begin{prop}\label{prop:age}
A class of finite relational structures $\mathcal C$ is closed under substructures, isomorphism, has the joint embedding property, and has countably many non-isomorphic members  
if and only if there exists a countable
structure $\Gamma$ whose age equals 
$\mathcal C$. 
\end{prop}

Proposition~\ref{prop:age} is the main motivation
why we exclusively work with classes of structures that are closed under substructures;
however, as demonstrated in a recent paper by
Zucker~\cite{Zucker14}, several Ramsey results 
and techniques can meaningfully be extended to isomorphism-closed classes that only satisfy the joint embedding property and amalgamation, but  that are not necessarily closed under substructures.

\begin{Def}[Ramsey class]
Let $\tau$ be an at most countable relational signature. 
A class of finite $\tau$-structures 
is called a \emph{Ramsey class} if it is closed under 
substructures, isomorphism, has countably many non-isomorphic members, the joint embedding, and the Ramsey property.
\end{Def}

%\begin{example}
%The classes that we have seen in Example~\ref{expl:all-structs-Ramsey}, or more generally, in Example~\ref{expl:irreducible-forbidden-Ramsey}, are Ramsey classes. 
%\end{example}

Examples of Ramsey classes will be presented below, in Example~\ref{expl:all-structs-Ramsey}, or more generally, in Example~\ref{expl:irreducible-forbidden-Ramsey}. 
The following can be shown by a simple compactness argument.

\begin{prop}\label{prop:compactness}
Let $\Gamma$ be a structure of age
$\mathcal C$. Then $\mathcal C$ is a Ramsey class if and only if 
for all $\bA,\bB \in \mathcal C$ and $r \in \mN$ we have that
$\Gamma \rightarrow (\bB)^\bA_r$.
\end{prop}
\begin{proof}
Let $\bA,\bB \in \mathcal C$, and $r \in \mN$ an integer. 
%When $\mathcal C$ is Ramsey, then there exists a finite $\bC \in \mathcal C$ such that  $\bC \rightarrow (\bB)^\bA_r$,
%and because $\bC$ is a substructure of $\Gamma$, we also have $\bC \rightarrow (\bB)^\bA_r$. 
%Conversely, suppose that $\Gamma \rightarrow (\bB)^\bA_r$.
When $k$ is the cardinality of ${\bB \choose \bA}$, then for any structure $\bC$ 
the fact that $\bC \rightarrow (\bB)^\bA_r$ can equivalently be expressed in terms of $r$-colourability of a certain $k$-uniform hypergraph, defined as follows. Let $G=(V;E)$ be the structure whose vertex set $V$ is ${\bC \choose \bA}$,
and where $(e_1,\dots,e_k) \in E$ if there exists an $f \in {\bC \choose \bB}$ such that $f \circ {\bB \choose \bA} = \{e_1,\dots,e_k\}$.
Let $H = ([r];E)$ be the structure where $E$ contains all tuples except for the tuples $(1,\dots,1),\dots,(r,\dots,r)$. 
Then $\bC \not\rightarrow (\bB)^\bA_r$ if and only if $G$ does not homomorphically map to $H$.
An easy and well-known compactness argument (see Lemma 3.1.5 in~\cite{Bodirsky-HDR-v4}) shows that 
this is the case if and only if some finite substructure of $G$ does not homomorphically map to $H$. 
Thus, $\Gamma \rightarrow (\bB)^\bA_r$ 
if and only if 
$\bC \rightarrow (\bB)^\bA_r$ for all finite substructures $\bC$ of $\Gamma$.
\end{proof}

\subsection{Ramsey degrees and rigidity}
\label{sect:rigidity}
Let $\mathcal C$ be a class of structures with the Ramsey property. 
In this section we will see that 
each structure in $\mathcal C$ must
be \emph{rigid}, that is, it has no automorphism other than the identity. 
%\begin{proof}
%Suppose that $B \in \mathcal C$ is such that 
%$d := |\Aut(B)| > 1$. 
%Let $C \in \mathcal C$ be such that $B$ is a %substructure of $C$. Let $<^C$ be a linear %ordering of $C$. Define $\chi \colon {C \choose %B} \to \Aut(B)$ as follows.
%For $e \in {C \choose B}$, define $\chi(e) = e$ 
%if 
%satisfies $C \to (B)^B_2$: 
%\end{proof}
%We show something slightly more general,
%using a notion that later
%becomes important again 
%in Section~\ref{sect:expansions}. 

\begin{defn}[Ramsey degrees]
Let $\mathcal C$ be a class of structures and let
$\bA \in \mathcal C$. We say that $\bA$ has \emph{Ramsey degree $k$ (in $\mathcal C$)} if $k \in \mN$ is least such that for any %superstructure 
$\bB \in \mathcal C$ %of $\bA$ 
and for any $r \in \mN$ there exists a $\bC \in \mathcal C$ such that for any $r$-colouring $\chi$ of ${\bC \choose \bA}$ there is an $f \in {\bC \choose \bB}$ such that $|\chi(f \circ {\bB \choose \bA})| \leq k$. 
\end{defn}

Hence, by definition, $\mathcal C$ has the
Ramsey property if every $\bA \in \mathcal C$
has Ramsey degree one. 

\begin{lemma}
Let $\mathcal C$ be a class of finite structures. 
Then for every $\bA \in \mathcal C$,
the Ramsey degree of $\bA$ in $\mathcal C$
is at least $|\Aut(\bA)|$. 
\end{lemma}
\begin{proof}
We have to show that for some $\bB \in \mathcal C$ and $r \in \mathbb N$,
every $\bC \in \mathcal C$ can be $r$-coloured such that for all $f \in {\bC \choose \bB}$ we have $|\chi(f \circ {\bB \choose \bA})| \geq |\Aut(\bA)|$. We choose $\bB := \bA$ and $r:=|\Aut(\bA)|$. 

Let $\bC \in \mathcal C$ be arbitrary. 
Define an equivalence relation $\sim$ on
${\bC \choose \bA}$ by setting $f \sim g$ if there
exists an $h \in \Aut(\bA)$ such that $f = g \circ h$. 
Let $f_1,\dots,f_t$ be a list of representatives
for the equivalence classes of $\sim$. 
Define $\chi \colon {\bC \choose \bA} \to \Aut(\bA)$ as follows. For $f \in {\bC \choose \bA}$, let $i$ be the unique $i$ 
such that $f_i \sim f$. Define $\chi(f) = h$ if 
$f = f_i \circ h$. Now let $e \in {\bC \choose \bA}$ be arbitrary. Then $|\chi(e \circ {\bA \choose \bA})| = |\Aut(\bA)|$. 
\end{proof}

\begin{cor}\label{cor:rigid}
Let $\mathcal C$ be a class with the Ramsey property. Then all $\bA$ in $\mathcal C$
are rigid.
\end{cor}

It follows that in particular the class of all finite graphs does \emph{not} have the Ramsey property. 
Frequently, 
a class without the Ramsey property
can be made Ramsey by expanding its members appropriately 
with a linear ordering (the expanded structures are
clearly rigid). 

\begin{example}\label{expl:all-structs-Ramsey}
Abramson and Harrington~\cite{AbramsonHarrington} and independently \Nesetril\ and R\"odl~\cite{NesetrilRoedlOrderedStructures}  
showed that for any relational signature $\tau$, the class $\cal C$ of all finite \emph{linearly ordered} $\tau$-structures has the Ramsey property. That is, the members of $\cal C$ are finite structures 
$\bA = (A; \preceq, R_1,R_2,\dots)$ for some fixed signature $\tau = \{\preceq,R_1,R_2,\dots\}$ 
where $\preceq$ denotes a linear order of $A$.

A shorter and simpler proof of this substantial result, based on the \emph{partite method}, can 
be found in~\cite{NesetrilRoedlPartite} and \cite{NesetrilSurvey} and will be presented in Section~\ref{sect:partite}. \qed
\end{example} 

For a class of finite $\tau$-structures $\mathcal N$, we write $\Forb({\mathcal N})$ for the class
of all finite $\tau$-structures that do not admit a homomorphism from any structure in $\mathcal N$.

\begin{example}\label{expl:irreducible-forbidden-Ramsey}
The classes from Example~\ref{expl:all-structs-Ramsey} have been further generalised by \Nesetril\ and R\"odl~\cite{NesetrilRoedlOrderedStructures} as follows.
Suppose that $\mathcal N$ is a (not necessarily finite) class of structures $\bF$
with finite 
relational signature $\tau$ such that 
for all elements $u,v$ of $\bF$ there is a tuple
in a relation $R^\bF$ for $R \in \tau$ that contains both $u$ and $v$. Such structures
have been called \emph{irreducible} in the Ramsey theory literature. 
%It can be readily verified that $\C $ is an amalgamation class. 
Then the class of all
expansions of the structures in $\C := \Forb(\mathcal N)$ by a linear order has the Ramsey property.
Again, there is a proof based on the partite method, which will be presented in Section~\ref{sect:partite}. 
This is indeed a generalization since we obtain the classes from Example~\ref{expl:all-structs-Ramsey} by taking ${\mathcal N} = \emptyset$. \qed
\end{example}

\subsection{The amalgamation property}
\label{sect:amalgamation}
The Ramsey classes we have seen so far 
will look familiar to model theorists.  
As mentioned in the introduction, the fact that all of the above Ramsey classes could be described as the age of a homogeneous structure is not a coincidence.

\begin{theorem}[\cite{RamseyClasses}]\label{thm:ramsey-amalgamation}
Let $\tau$ be a relational signature, and let $\mathcal C$ be a class of 
finite $\tau$-structures that 
is closed under 
isomorphism, 
and has the joint embedding property. 
If $\mathcal C$ has the Ramsey property,
then it also has the amalgamation property.
\end{theorem}
\begin{proof}
Let $\bA,\bB_1,\bB_2$ be members of $\mathcal C$
such that there are embeddings $e_i \in {\bB_i \choose \bA}$ for $i = 1$ and $i = 2$. Since $\cal C$ has the joint embedding property, 
there exists a structure $\bC \in \cal C$
with embeddings $f_1,f_2$ of $\bB_1$ and $\bB_2$ into $\bC$.
%If the restriction of $e_1$ to $A$ equals the
%restriction of $e_2$ to $A$, 
If $f_1 \circ e_1 = f_2 \circ e_2$, 
then $\bC$ shows that $\bB_1$ and $\bB_2$ amalgamate over $\bA$, so
assume otherwise. 

Let $\bD \in \mathcal C$ be such that $\bD \rightarrow (\bC)^{\bA}_2$. 
Define a colouring 
$\chi \colon {\bD \choose \bA} \rightarrow [2]$ as follows. 
For $g \in {\bD \choose \bA}$, 
let $\chi(g)=1$ if there
 is a $t \in {\bD \choose \bC}$ 
such that $g = t \circ f_1 \circ e_1$, and $\chi(g)=0$ otherwise. 
Since $\bD \rightarrow (\bC)^{\bA}_2$, there exists a $t_0 \in {\bD \choose \bC}$
such that $|\chi(t_0 \circ {\bC \choose \bA})| = 1$. 
%For $i \in \{1,2\}$, 
%let $e_i'$ be the restriction of $e_i$ to $A$.
Note that $\chi(t_0 \circ f_1 \circ e_1) = 1$ by the definition of $\chi$. It follows that $\chi(t_0 \circ h) = 1$ 
for all $h \in {\bC \choose \bA}$. 
In particular $\chi(t_0 \circ f_2 \circ e_2) = 1$, because
$f_2 \circ e_2 \in {\bC \choose \bA}$. 
Thus, by the definition of $\chi$, 
there exists a $t_1 \in {\bD \choose \bC}$
such that $t_1 \circ f_1 \circ e_1 = t_0 \circ f_2 \circ e_2$ (here we use that the structure $\bA$ must be rigid, by Corollary~\ref{cor:rigid}). 
This shows that $\bD$ together with the embeddings 
$t_1 \circ f_1 \colon \bB_1 \rightarrow \bD$ 
and $t_0 \circ f_2 \colon \bB_2 \rightarrow \bD$ 
is an amalgam of $\bB_1$ and $\bB_2$ over $\bA$.
\end{proof}

\begin{Def}[Amalgamation Class]\label{def:amalgamation-prop}
An isomorphism-closed class of finite
structures with an at most countable relational
signature 
that contains at most countably many non-isomorphic structures,
has the amalgamation property (defined in the introduction), and that is closed under
taking induced substructures,
is called an \emph{amalgamation class}.
\end{Def}

\begin{theorem}[\Fresse~\cite{OriginalFraisse,Fraisse}; see~\cite{Hodges}]\label{thm:fraisse}
Let $\tau$ be a countable relational signature and let 
$\cal C$ be an amalgamation class of $\tau$-structures. 
Then there is a homogeneous and at 
most countable $\tau$-structure $\bC$ whose age equals $\cal C$.
The structure $\bC$ is unique up to isomorphism, and called
the \emph{\Fresse\ limit} of $\cal C$.
\end{theorem}

\begin{example}
The \Fresse\ limit of the class of all finite linear orders is isomorphic to $({\mathbb Q};<)$, the order of the rationals. 
The \Fresse\ limit of the class of all graphs
is the so-called random graph (or Rado graph); see e.g.~\cite{RandomRevisitedCameron}. 
\end{example}

We also have the following converse of 
Theorem~\ref{thm:fraisse}. 

\begin{theorem}[\Fresse; see~\cite{Hodges}]\label{thm:fraisse2}
Let $\Gamma$ be a homogeneous
relational structure. Then the age of
$\Gamma$ is an amalgamation class. 
\end{theorem}

As we have seen, there is a close connection between
amalgamation classes and homogeneous structures, and
we therefore make the following definition. 

\begin{Def}[Ramsey structure]\label{def:ramsey}
A homogeneous 
structure $\Gamma$ is called \emph{Ramsey}
if the age of $\Gamma$ has the Ramsey property.
\end{Def}

\subsection{Counterexamples}
\label{sect:counter}
We have so far seen two important necessary
conditions for a class $\mathcal C$ to
be a Ramsey class: rigidity of the members of $\mathcal C$ (Corollary~\ref{cor:rigid})
and amalgamation (Theorem~\ref{thm:ramsey-amalgamation}). 
As we will see in the examples in this section,
these conditions are not sufficient for being 
Ramsey.

%\noindent A C-relation is called \emph{dense} if it satisfies
%\begin{enumerate}
%\item [C7] $C(a; b,c) \rightarrow \exists e \, (C(e;b,c) \wedge C(a; b,e))$.
%\end{enumerate}

\begin{example}\label{expl:equivalences}
Let $\mathcal C$ be the class of all finite $\{E,<\}$-structures where $E$ denotes an equivalence relation and $<$ denotes a linear order.
It is easy to verify that $\mathcal C$ has the
amalgamation property. Moreover, all automorphisms of structures in $\mathcal C$ have to preserve $<$ and hence must be the identity. But $\mathcal C$ does not have the Ramsey property: 
let $\bA$ be the structure with domain $\{u,v\}$ such that $<^\bA \, = \{(u,v)\}$,
and such that $u$ and $v$ are not $E$-equivalent. Let $\bB$ be the structure with domain $\{a,b,c,d\}$
such that $b <^\bB c <^\bB a <^\bB d$ and such that $\{a,b\}$ and $\{c,d\}$ are the equivalence classes of $E^\bB$.  There are four copies of $\bA$ in $\bB$. 

Suppose for contradiction that 
there is $\bC \in \mathcal C$ 
such that $\bC \rightarrow (\bB)^\bA_2$. 
Let $\prec$ be a \emph{convex} linear ordering of the elements
of $C$, that is, a linear ordering such that
$E(x,z)$ and $x < y < z$ implies that $E(x,y)$
and $E(y,z)$. 
Let $g \in {\bC \choose \bA}$. Define $\chi(g) = 1$ if $g(u) \prec g(v)$, 
and $\chi(g) = 2$ otherwise. 
Note that there are only two convex linear
orderings of $\bB$, and that $|\chi(f \circ {\bB \choose \bA})| = 2$ for all $f \in {\bC \choose \bB}$. 
\qed
\end{example}

However, the class of all equivalence relations
with a \emph{convex} linear order is Ramsey; 
see~\cite{Topo-Dynamics}. % (this can also be shown with the general tools presented in Section...)
Moreover, 
as we will see in Example~\ref{expl:E-with-two-orders} in Section~\ref{sect:superimposing}, 
the \Fresse\ limit of the class $\mathcal C$ from Example~\ref{expl:equivalences} can be expanded by a convex linear order $\prec$ so that the resulting structure is homogeneous and Ramsey.

\ignore{
\begin{example}\label{expl:s2}
Consider the class of all partially ordered
finite sets $(S;\leq)$ such that there exists a
finite rooted tree whose vertices contain $S$ such
that $u \leq v$ for $u,v \in S$ if and only if 
$v$ lies on the path from $u$ to the root. 
This class does not have amalgamation, 
but the class of all finite structure $(S;\leq,C)$ 
does, where $C$ is the relation that contains 
all triples $(a,b,c)$ of element from $S$ such that
%A partially ordered set is called \emph{semilinear} if 
\end{example}
}

\begin{figure}
\begin{center}
  	\includegraphics[scale=0.4]{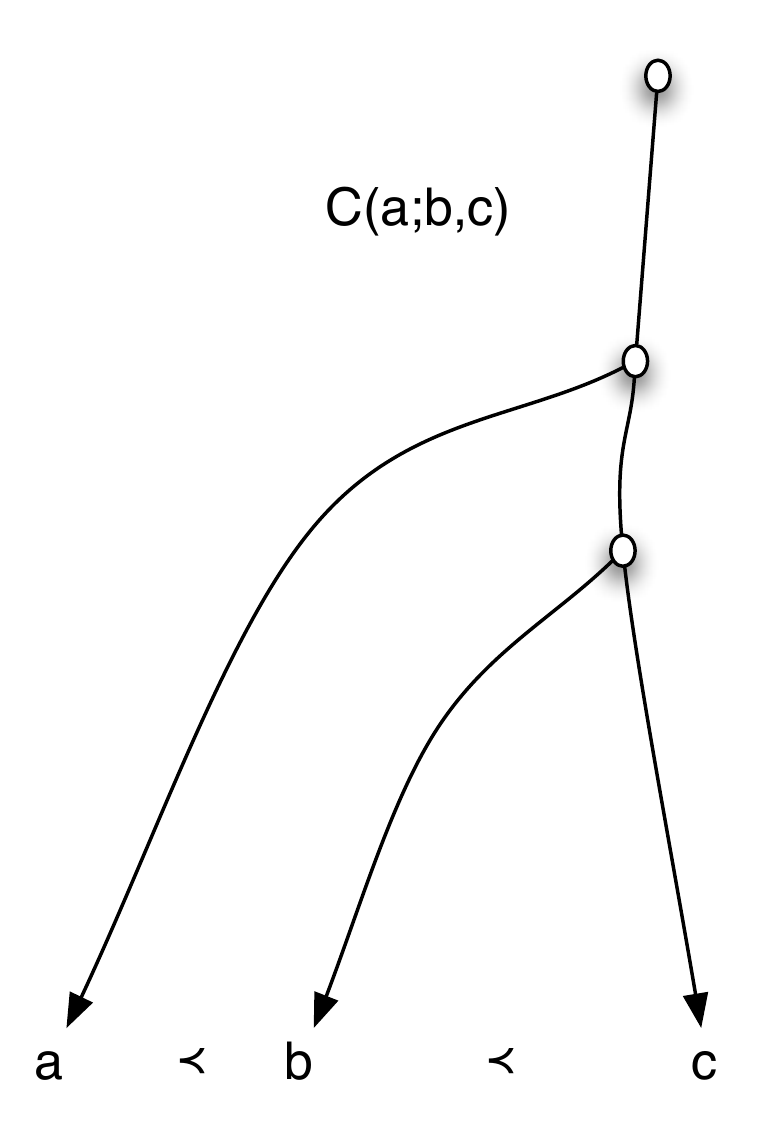} 
	\end{center}
\caption{Illustration of a convexly ordered C-relation.}
\label{fig:C}
\end{figure}

\begin{example}\label{expl:c-rel}
The class of finite trees is not closed under taking substructures. If we close it under substructures, we
obtain the class of all finite forests, a class which does not have the amalgamation property. 
The solution for a proper model-theoretic treatment
of trees and forests 
is to use the concept of $C$-relations. 

Formally, a ternary relation $C$ is said to be a \emph{C-relation}\footnote{Terminology of Adeleke and Neumann~\cite{AdelekeNeumann}.} on a set $L$ if for all $a,b,c,d \in L$ the following conditions hold:
\begin{enumerate}
\item [C1] $C(a; b,c) \rightarrow C(a; c,b)$; 
\item [C2] $C(a; b,c) \rightarrow \neg C(b;a,c)$; 
\item [C3] $C(a; b,c) \rightarrow C(a; d,c) \vee C(d; b,c)$; 
\item [C4] $a \neq b \rightarrow C(a;b,b)$.
\end{enumerate}
A $C$-relation on a set $L$ is called \emph{binary branching} if for all pairwise distinct $a,b,c \in L$ we have $C(a;b,c)$ or $C(b;a,c)$
or $C(c;a,b)$. 

The intuition here is that the elements of $L$ denote the leaves of a rooted binary tree, and $C(a;b,c)$ holds if in the tree, the shortest path from $b$ to $c$ does not intersect the shortest path from $a$ to the root; see Figure~\ref{fig:C}. 
For finite $L$, %(or when we require this for all finite subsets of $L$)
this property is actually equivalent to the axiomatic definition above~\cite{AdelekeNeumann}. 

The class of structures $(L;C)$ where $L$ is a finite set and $C$ is a binary branching 
C-relation on $L$ 
is of course not a Ramsey class, since 
$(L;C)$ has nontrivial automorphisms, 
unless $|L| = 1$. 
The same argument does not work for the class $\mathcal C$ 
of all structures $(L;C,<)$ where $L$ is finite set, $C$ is a binary branching C-relation on $L$, and $<$ is a linear ordering of $L$. In fact, $\mathcal C$ is an amalgamation class (a well-known fact; for a proof, see~\cite{Bodirsky-HDR-v4}),
%Section~\ref{sect:inductive}), 
but not a Ramsey class.
To see how the Ramsey property fails, 
consider the structure $\bB \in \mathcal C$
with domain $\{a,b,c,d\}$ where $a<c<b<d$ such that
 $C(a;c,d),C(b;c,d),C(d;a,b),C(c;a,b)$, 
and the structure $\bA \in \mathcal C$ with domain
$\{u,v\}$ where $u < v$. Now let $\bC \in \mathcal C$ 
be arbitrary. Let $\prec$ be a \emph{convex}
ordering of $\bC$, that is, a linear ordering
such that for all $u,v,w \in L$, 
if $C(u;v,w)$ and $v \prec w$, then either $u \prec v \prec w$
or $v \prec w \prec u$. Define $\chi \colon {\bC \choose \bA} \to [2]$ as follows. For $g \in {\bC \choose \bA}$
define $\chi(g) = 1$ if $g(u) \prec g(v)$,
and $\chi(g) = 2$ otherwise. Note that for 
every convex ordering $\prec$ of $B$ there
exists an $e_1 \in {\bB \choose \bA}$ such that
$e_1(u) \prec e_1(v)$, and an 
$e_2 \in {\bB \choose \bA}$ such that
$e_2(v) \prec e_2(u)$. Hence, 
for every $f \in {\bC \choose \bB}$ we have
$|\chi(f \circ {\bB \choose \bA})| = 2$. 
\end{example}

Again, the class of all convexly ordered binary branching $C$-relations over a finite set is an amalgamation class (Theorem~\ref{thm:tramsey}). 
Moreover, by the results from Section~\ref{sect:superimposing}, 
the \Fresse\ limit of the class $\mathcal C$ from Example~\ref{expl:c-rel} can be expanded by a convex linear order so that the resulting structure is homogeneous and Ramsey; see Example~\ref{expl:C-with-two-orders}.
% in Section~\ref{sect:inductive}. 

\subsection{Automorphism groups}
\label{sect:automorphisms}
Let $f \colon D \to D$ be a function 
and $t \in D^m$ a tuple. 
Then $f(t)$ denotes the tuple $(f(t_1),\dots,f(t_m))$. 
We say that a relation $R \subseteq D^m$ 
is \emph{preserved} by a function $f \colon D \to D$ if $f(t) \in R$ for all $t \in R$. 
An automorphism of a structure
$\Gamma$ with domain $D$ is a permutation
$\alpha$ such that both $\alpha$ and $\alpha^{-1}$ preserve all relations (and if the signature contains constant symbols,
$\alpha$ must fix the constants).

%Proposition~\ref{prop:group} below
%shows that whether a homogeneous
%structure $\Gamma$ 
%is Ramsey only depends 
%on the automorphism group\footnote{We view the automorphism group as a permutation group; in fact, whether a homogeneous structure $\Gamma$ is Ramsey or not only depends on the automorphism group of $\Gamma$, viewed as a \emph{topological group}; see~\cite{Topo-Dynamics}.} 
%$\Aut(\Gamma)$ of $\Gamma$.
The equivalent formulation of the Ramsey property in Proposition~\ref{prop:clean} 
will be useful later, for instance 
to prove that
for every homogeneous Ramsey structure 
$\Gamma$ there exists a linear
order on the domain of $\Gamma$ that
is preserved by all automorphisms of $\Gamma$. 

\begin{prop}\label{prop:clean}
Let $\Gamma$ be a homogeneous structure. 
Then the following are equivalent.
\begin{enumerate}
\item $\Gamma$ is Ramsey.
\item For every finite substructure $\bB$ of $\Gamma$ and $r \in \mN$ there exists a finite 
substructure $\bC$ of $\Gamma$
such that for all substructures $\bA_1,\dots,\bA_\ell$ of $\bB$ and all
$\chi_i \colon {\bC \choose \bA_i} \to [r]$
 there exists an $e \in {\bC \choose \bB}$ 
such that $|\chi_i(e \circ {\bB \choose \bA_i})| = 1$ for all $i \in [\ell]$. 
%for all $e \in {\bB \choose \bA}$
%and $\alpha \in \Aut(\Gamma)$ we have $\chi(\beta \circ e) = \chi(\alpha \circ \beta \circ e)$. 
%\item For every finite substructure $\bB$ of $\Gamma$ and $r \in \mN$ and $\chi \colon {D \choose S} \to [r]$ there exist $\beta \in \Aut(\Gamma)$ 
%such that for all $e \in {\beta(M) \choose S}$ and $\alpha \in \Aut(\Gamma)$ we have 
%that $\chi(\alpha \circ e) = \chi(e)$. 
\end{enumerate}
\end{prop}
\begin{proof}
$(1) \Rightarrow (2)$. 
We only show the forward implication, the backward implication being trivial. 
Let $\bB$ be a finite substructure of $\Gamma$ and $r \in \mN$. 
%, and $\chi \colon {D \choose S} \to [r]$ be given. 
Let $\bA_1,\dots,\bA_\ell$ be an enumeration 
of the substructures of $\bB$.
We are going to construct a sequence of
structures $\bC_1,\dots,\bC_\ell$. Since $\Gamma$ is Ramsey, there exists
a substructure $\bC_1$ of $\Gamma$ 
such that $\bC_1 \to (\bB)^{\bA_1}_{r}$. 
Inductively, for $i \in \{2,\dots,\ell\}$ there exists a substructure $\bC_i$ of $\Gamma$
such that $\bC_i \to (\bC_{i-1})^{\bA_i}_{r}$. 
Define $\bC := \bC_\ell$. 

%Now, let $\bA$ be an arbitrary substructure of $\bB$, and $\chi \colon {\bC_l \choose \bA} \to [r]$ be arbitrary. Then $\bA = \bA_i$ for some $i \in [l]$, and since $\bC_i \to (\bC_{i-1})^{\bA_i}_{r}$,
%there exists an $e_i \in {\bC_i \choose \bC_{i-1}}$ with $|\chi(e \circ {\bC_{i-1} \choose \bA_i})| \leq 1$. For all $j \in [l]$, $\bC_{j-1}$ is a substructure of $\bC_j$, and in particular $\bC_i$ is a substructure of $\bC_l$ and $\bB = \bC_0$ is a substructure of $\bC_{i-1}$. 

For all $i \in [\ell]$, let $\chi_i \colon {\bC \choose \bA_i} \to [r]$ be arbitrary. Since $\bC_\ell \to (\bC_{\ell-1})^{\bA_\ell}_r$, there exists an $e_\ell \in {\bC_\ell \choose \bC_{\ell-1}}$ with $|\chi(e_\ell \circ {\bC_{\ell-1} \choose \bA_\ell})| \leq 1$. 
Inductively, suppose we have already
defined $e_i \in {\bC_\ell \choose \bC_{i-1}}$ for an $i \in \{2,\dots,\ell\}$ such that
for all $j \in \{i,\dots,\ell\}$ we have 
$|\chi(e_i \circ {\bC_{i-1} \choose \bA_j})| \leq 1$.
Then there exists an $e_{i-1} \in {\bC_{i-1} \choose \bC_{i-2}}$ such that $|\chi(e_{i-1} \circ {\bC_{i-2} \choose \bA_{i-1}})| \leq 1$. Hence, for all $j \in \{i-1,\dots,\ell\}$ we have $|\chi(e_{i-1} \circ {\bC_{i-2} \choose \bA_j})| \leq 1$. Then the map 
$e_1 \in {\bC_k \choose \bB}$ has 
%By the homogeneity of $\Gamma$, the map
%$e_1 \in {\bC_k \choose \bB}$ can be extended to an automorphism $\beta$ of $\Gamma$, and this automorphism has
the desired properties from the statement of the proposition.
%$(2) \Rightarrow (1)$: trivial. 
%$(2) \Rightarrow (3)$. Trivial.
%$(3) \Rightarrow (1)$. 
%Let $A,B$ be finite substructures of $\Gamma$, and $r \in \mN$. Let $k := |A|$. 
%A straightforward consequence of Proposition~\ref{prop:compactness}. 
%it suffices to show that $\Gamma \to (B)^A_r$. 
%So let $\chi \colon {\Gamma \choose A} \to [r]$ be arbitrary. 
%Define $\xi \colon {D \choose [k]} \to [r]$ as follows. If the image of 
%$e \in {D \choose [k]}$ induces a copy of $A$ in $\Gamma$, then $\xi(e) := \chi(e)$.
%Otherwise, define $\xi(e)$ arbitrarily. 
%Then the assumptions imply that there exists
%a $\beta \in \Aut(\Gamma)$ such that 
%for all $K_1,K_2 \subseteq \beta(B)$ with
%$|K_1| = |K_2| = k$ and $\alpha \in \Aut(\Gamma)$, if $\alpha(K_1) = K_2$ then 
%$\xi(K_1) = \xi(K_2)$. We claim that 
%the structure $C$ induced in $\Gamma$ by
%$\beta(B)$ has the property that $C \to (B)^A_r$. To see this, let $e_1,e_2 \in {C \choose A}$ be arbitrary. Let $K_1$ be the image of $e_1$
%and $K_2$ be the image of $e_2$. 
%By homogeneity of $\Gamma$ 
%there exists an $\alpha \in \Aut(\Gamma)$ such that $\alpha(K_1) = \alpha(K_2)$, and 
%hence $\xi(K_1) = \xi(K_2)$. Therefore,
%$\chi(e_1) = \chi(e_2)$, which is what we had to show. 
\end{proof}

\begin{prop}\label{prop:ordering-invariant}
Let $\Gamma$ be a homogeneous
Ramsey structure with domain $D$. Then there exists a linear order on $D$ that is preserved by
all automorphisms of $\Gamma$. 
\end{prop}
\begin{proof}
Let $d_1,d_2,\dots$ be an enumeration of $D$,
and let $<$ be the linear order on $D$ given by
this enumeration, that is, $d_i < d_j$ if and only if $i < j$. 
Let $\mathcal T$ be a tree whose vertices on level $n$ 
are linear orders $\prec$ of $D_n := \{d_1,\dots,d_n\}$
with the property that for all $a,b \in D_n$ and 
$\alpha \in \Aut(\Gamma)$ such that $\alpha(a),\alpha(b) \in D_n$, 
we have that $a \prec b$ if and only if $\alpha(a) \prec \alpha(b)$.
Note that when a linear order satisfies this
condition, then also restrictions of the linear order to subsets satisfy this condition. 
Adjacency in $\mathcal T$ is defined by restriction. Clearly, $\mathcal T$ is finitely branching. We will show that $\mathcal T$ has vertices on each
level. By K\"onig's lemma, there is an infinite path in $\mathcal T$,
which defines a linear ordering on $D$
that is preserved by $\Aut(\Gamma)$. 

To show that there is a linear order $\prec$ on
$D_n$ that satisfies the condition, 
let $\bB$ be the structure induced by $D_n$ in $\Gamma$, and let $\bA_1,\dots,\bA_\ell$ list the substructures of $\Gamma$ that are induced by the two-element subsets of $D_n$.
By Proposition~\ref{prop:clean}, there
exists a finite substructure $\bC$ of
$\Gamma$ such that for all $\chi_i \colon {\bC \choose \bA_i} \to [r]$ there exists an
$e \in {\bC \choose \bB}$ such that 
$|\chi_i(e \circ {\bB \choose \bA_i})| = 1$ for all $i  \in [l]$. Let
 $\chi_i \colon {\bC \choose \bA_i} \to [2]$ 
 be defined as follows. 
For $e \in {\bC \choose \bA_i}$, 
%and $A_i = \{d_p,d_q\}$ with $p < q$, 
%let $u,v \in \mN$ be such that $e(d_p) = d_u$
%and $e(y) = d_v$. 
we define $\chi_i(e) = 1$ if $e$ preserves $<$,
and $\chi_i(e) = 2$ otherwise. 
%$u<v$ 
%and $\chi_i(e) = 2$ otherwise. 
By the property of $\bC$, there is
an $e \in {\bC \choose \bB}$ be such that 
$|\chi_i(e \circ {\bB \choose \bA_i})| = 1$ for all $i  \in [l]$. 

Let $\prec$ be the linear order on $D_n$ given by $a \prec b$ if $e(a) < e(b)$. 
%$e(a) = d_u$,
%$e(b) = d_v$, and $u<v$. 
Suppose now
that $a,b \in D_n$ and $\alpha \in \Aut(\Gamma)$ such that $\alpha(a),\alpha(b) \in D_n$.
Let $i \in [l]$ be such that $\{a,b\}$ induce $\bA_i$ in $\Gamma$. Let $f_1$ be the identity on $\{a,b\}$, and let $f_2$ be the restriction of
$\alpha$ to $\{a,b\}$; then 
$f_1,f_2 \in {\bB \choose \bA_i}$, and
$\chi_i(e \circ f_1) = \chi_i(e \circ f_2)$. 
By the definition of $\chi_i$, we have
that $e(a) < e(b)$ if and only if $e(\alpha(a)) < e(\alpha(b))$. By the definition of $\prec$
we obtain that $a \prec b$ if and only if $\alpha(a) \prec \alpha(b)$. 
\end{proof}

\subsection{Countably categorical structures}
\label{sect:omega-cat}
In this subsection we present a generalization
of the class of all homogeneous structures with a finite relational signature that still satisfies
a certain finiteness condition, namely
the class of all countable \emph{$\omega$-categorical} structures. 

\begin{Def}
A countable structure is said to be \emph{$\omega$-categorical} if all
countable structures that satisfy 
the same first-order sentences as $\Gamma$
are isomorphic to $\Gamma$. 
\end{Def}

Theorem~\ref{thm:ryll} below explains why
$\omega$-categoricity can be seen as a finiteness condition. 
It will be easy to see from Theorem~\ref{thm:ryll} that all structures that are homogeneous in a finite relational signature are $\omega$-categorical. But we first show
an example where $\omega$-categoricity can be seen 
directly. 

\begin{example}
All countably infinite vector spaces $\mathfrak V$ over a fixed finite field $\mathbb F$ are isomorphic. Since the isomorphism type of $\mathbb F$,
the axioms of vector spaces, and having infinite dimension 
can be expressed by first-order sentences it follows that $\mathfrak V$
is $\omega$-categorical. 
These structures are homogeneous; however,
their signature is not relational. 
The relational structure with the same domain
that contains all relations that are first-order definable over $\mathfrak V$ is homogeneous, too (this follows from Theorem~\ref{thm:ryll} below). It is easy to
see that all relational 
structures obtained from those
examples by dropping all but finitely many relations, but have the same automorphism group as $\mathfrak V$, are \emph{not} homogeneous. 
For example, the structure that just contains
the ternary relation defined by $x = y+z$ 
has the same automorphism group as $\mathfrak V$, but is
\emph{not} homogeneous. 
The Ramsey properties of those examples are 
beyond the scope of this survey, but are discussed in~\cite{Topo-Dynamics}. 
\end{example}

The following theorem of Engeler, Ryll-Nardzewski, and Svenonius shows that whether a structure is $\omega$-categorical can be seen from the automorphism group $\Aut(\Gamma)$ of $\Gamma$ (as a permutation group). 

\begin{theorem}[see e.g.~\cite{Hodges}]
\label{thm:ryll}
Let $\Gamma$ be a countably infinite structure with a countably infinite signature. Then the following are equivalent. 
\begin{enumerate}
\item $\Gamma$ is $\omega$-categorical; 
\item $\Aut(\Gamma)$ is \emph{oligomorphic}, that
is, for all $n \geq 1$, the componentwise action
of $\Aut(\Gamma)$ on $n$-tuples from $\Gamma$ has finitely many orbits;
\item all orbits of $n$-tuples in $\Gamma$
are first-order definable in $\Gamma$;
\item all relations preserved by $\Aut(\Gamma)$ are first-order definable in $\Gamma$. 
\end{enumerate}
\end{theorem}

The following is a direct consequence
of Proposition~\ref{prop:ordering-invariant}
and Theorem~\ref{thm:ryll}. 

\begin{cor}\label{cor:omega-cat-order}
Let $\Gamma$ be an $\omega$-categorical
Ramsey structure. Then there is a
linear order with a first-order definition in
$\Gamma$. 
\end{cor}
\begin{proof}
Let $\Gamma^*$ be the homogeneous expansion of $\Gamma$ by all first-order definable relations.
By Proposition~\ref{prop:ordering-invariant}, 
there exists a linear ordering of the domain of $\Gamma^*$ which is preserved by all automorphisms of $\Gamma$. By $\omega$-categoricity of $\Gamma$ and $\Gamma^*$, Theorem~\ref{thm:ryll}, this linear order is first-order definable in $\Gamma^*$. Since
all first-order definable relations of $\Gamma^*$ are first-order definable in $\Gamma$, they are present in the signature of $\Gamma^*$, and the statement follows. 
\end{proof}

%Homogeneous structures with a finite relational signature have an important model-theoretic property, the are \emph{$\omega$-categorical}.
%Several general procedures to derive new Ramsey
%classes from existing Ramsey classes can naturally
%be phrased in terms of $\omega$-categorical structures. Moreover, there is an important link between $\omega$-categorical structures and permutation group theory, which
%is useful in many respects; we therefore give
%a brief introduction. 
Theorem~\ref{thm:ryll} implies 
that when $\Gamma$ is $\omega$-categorical, then the expansion $\Gamma'$ 
of $\Gamma$ by all first-order definable relations is homogeneous. We therefore make
the following definition. 

\begin{Def}\label{def:omega-cat-ramsey}
An $\omega$-categorical structure $\Gamma$
is called \emph{Ramsey} if the expansion
of $\Gamma$ by all relations with a first-order definition in $\Gamma$ is Ramsey (as a homogeneous structure).
\end{Def}

This definition is compatible with Definition~\ref{def:ramsey}, since
expansions by first-order definable relations
do not change the automorphism group, and since the Ramsey property only depends on the automorphism group, as reflected in the next proposition. 
For subsets $S$ and $M$ of the domain
of an $\omega$-categorical structure $\Gamma$, 
we write ${M \choose S}$ for the set
of all maps from $S$ to $M$ that can be extended to an automorphism of $\Gamma$. 
The following is immediate from the definitions, Theorem~\ref{thm:ryll}, 
and Proposition~\ref{prop:compactness}. 

\begin{prop}\label{prop:group}
Let $\Gamma$ be an $\omega$-categorical structure with domain $D$. 
Then the following are equivalent. 
\begin{enumerate}
\item $\Gamma$ is Ramsey;
\item For all $r \in \mN$ and finite $M \subset D$ and $S \subset M$ there exists a finite $L \subseteq D$ such that for every map $\chi$ from ${L \choose S}$ to $[r]$ there exists $f \in {L \choose M}$ such that $|\chi(f \circ {M \choose S})| = 1$. 
\item For all $r \in \mN$ and finite $M \subset D$ and $S \subset M$ and every map $\chi$ from
${D \choose S}$ to $[r]$ there exists $f \in {D \choose M}$ such that $|\chi(f \circ {M \choose S})| = 1$. 
\end{enumerate}
\end{prop}
% TODO: proof!

%A structure $\Gamma$ has \emph{quantifier elimination} if every first-order formula 
%is over $\Gamma$ equivalent to a quantifier-free formula. 
%\begin{prop}\label{prop:qe}
%Let $\Gamma$ be $\omega$-categorical.
%Then $\Gamma$ is homogeneous if
%and only if $\Gamma$ has quantifier elimination. 
%\end{prop}

\section{New Ramsey classes from old}
\label{sect:new-old}
% FOCUSSING ON omega-cat HERE MIGHT
% SCARE THE READERS. This is why
% I have reformulated the following:
%The class of $\omega$-categorical structures
%is very robust with respect to basic model-theoretic constructions. When
%$\Gamma$ is $\omega$-categorical,
%then the following structures obtained from
%$\Gamma$ will also be $\omega$-categorical:
%\begin{itemize}
%\item reducts of $\Gamma$, and, more generally, structures with a first-order interpretation in $\Gamma$;
%\item expansions of $\Gamma$ by finitely many constants;
%\item the model companion of $\Gamma$;
%\item the model-complete core of $\Gamma$;
%\item under certain additional assumptions, $\Gamma$ can be \emph{superimposed} with another $\omega$-categorical structure, to obtain new $\omega$-categorical structures (precise definitions in Section~\ref{sect:superimposing}). 
%\end{itemize}
%The structures that we thus obtain will be again $\omega$-categorical. In this section we will see that if the original structures have good Ramsey properties, then the new structures also do. 

The class of $\omega$-categorical Ramsey structures is remarkably robust with respect to basic
model-theoretic constructions. We will consider the following 
model-theoretic constructions to obtain
new structures from given structures $\Gamma,\Gamma_1,\Gamma_2$: 
\begin{itemize}
\item disjoint unions and products of $\Gamma_1$ and $\Gamma_2$;
\item structures with a first-order interpretation in $\Gamma$;
\item expansions of $\Gamma$ by finitely many constants;
\item the model companion of $\Gamma$;
\item the model-complete core of $\Gamma$;
\item superpositions of $\Gamma_1$ and $\Gamma_2$.
\end{itemize}
If the structures $\Gamma,\Gamma_1,\Gamma_2$ we started from
are $\omega$-categorical (or homogeneous in a finite relational signature), the structure we thus obtain will be again $\omega$-categorical (or homogeneous in a finite relational signature). 
In this section we will see
that if the original structures have good Ramsey properties, then the new structures also do. 

\subsection{Disjoint unions}
One of the simplest operations on structures 
is the formation of disjoint unions:
when $\bA_1$ and $\bA_2$ are structures with
the same relational signature $\tau$ and
disjoint domains, then the disjoint union
of $\bA_1$ and $\bA_2$ is the structure $\bB$ with domain $B := A_1 \cup A_2$ where for each $R \in \tau$ we set $R^{\bB} := R^{\bA_1} \cup R^{\bA_2}$. 
The disjoint union of two $\omega$-categorical 
structures is always $\omega$-categorical.
The disjoint union of two 
homogeneous structures $\Gamma_1$ and $\Gamma_2$ might not be homogeneous; but it clearly becomes homogeneous when we add an additional new
unary predicate $P$ to the disjoint union which precisely contains the vertices from $\Gamma_1$. We denote the resulting structure 
by $\Gamma_1 \uplus_P \Gamma_2$. 
The transfer of the Ramsey property
is a triviality in this case. 
%The proof of the following lemma is straightforward and therefore omitted. 

\begin{lemma}\label{lem:disj-union}
Let $\Gamma_1$ and $\Gamma_2$ be 
$\omega$-categorical Ramsey structures. 
Then $\Gamma := \Gamma_1 \uplus_P \Gamma_2$ is an $\omega$-categorical Ramsey structure, too.
If $\Gamma_1$ and $\Gamma_2$ are homogeneous with finite relational signature, then so is $\Gamma$. 
\end{lemma}

While this lemma looks innocent, it still has interesting applications in combination
with the other constructions that we present; see Example~\ref{ex:s2}.

\subsection{Products}
When $G_1$ and $G_2$ are permutation groups acting on the sets $D_1$ and $D_2$, respectively, then the direct product $G_1 \times G_2$ of $G_1$ and $G_2$ naturally acts on $D_1 \times D_2$:
the element $(g_1,g_2)$ of $G_1 \times G_2$ maps $(x_1,x_2)$ to $(g_1(x_1),g_2(x_2))$. 
When $G_1$ and $G_2$ are the automorphism 
groups of relational structures $\Gamma_1$ and $\Gamma_2$, then the following definition 
yields a structure whose automorphism
group is precisely $G_1 \times G_2$. (The direct product $\Gamma_1 \times \Gamma_2$ does not have this property.) 

\begin{Def}[Full Product]
Let $\Gamma_1,\dots,\Gamma_d$ be structures with domains $D_1,\dots,D_d$ and pairwise disjoint signatures $\tau_1,\dots,\tau_d$. Then the \emph{full product structure} $\Gamma_1 \boxtimes \cdots \boxtimes \Gamma_d$ is the structure with domain $D_1 \times \cdots \times D_d$ that contains for every $i \leq d$ and $m$-ary $R \in (\tau_i \cup \{=\})$ the relation
defined by $\big \{ ((x_1^1,\dots,x_1^d),\dots,(x_m^1,\dots,x_m^d)) : (x_1^i,\dots,x_m^i) \in R^{\Gamma_i}  \big\}$. 
\end{Def}

The following proposition is known
as the \emph{product Ramsey theorem}
to combinatorists.

\begin{prop}\label{prop:product-ramsey}
Let $\Gamma_1,\dots,\Gamma_d$ be $\omega$-categorical Ramsey structures with pairwise disjoint signatures.
Then $\Gamma := \Gamma_1 \boxtimes \cdots \boxtimes \Gamma_d$ is $\omega$-categorical and Ramsey. If $\Gamma_1,\dots,\Gamma_d$ are homogeneous with finite relational signature, then so is $\Gamma$. 
\end{prop}
\begin{proof}
It follows from Theorem~\ref{thm:ryll} that
if $\Gamma_1,\dots,\Gamma_d$ 
are $\omega$-categorical,
then $\Gamma_1 \boxtimes \cdots \boxtimes \Gamma_d$ is $\omega$-categorical. 

For the homogeneity of $\Gamma_1 \boxtimes \Gamma_2$, 
let $u_1 := (u^1_1,\dots,u^d_1)$, \dots, $u_m := (u^1_m,\dots,u^d_m)$ and $v_1 := (v^1_1,\dots,v^d_1)$, \dots, $v_m := (v^1_m,\dots,v^d_m)$ be elements of
$\Gamma$ such that the map $a$ 
that sends 
$(u_1,\dots,u_m)$ to $(v_1,\dots,v_m)$ is an isomorphism between
substructures of $\Gamma$. 
For $i \leq d$, define $a_i$ as the
map that sends $u^i_j$ to $v^i_j$ for all $j \leq m$; this is well-defined since $a$ preserves
the relation $\{(x^1,\dots,x^d,y^1,\dots,y^d) : x^i = y^i\}$. By homogeneity of $\Gamma_i$,
there exists an extension
$\alpha_i$ of $a_i$ to an automorphism of $\Gamma_i$.
Then the map $\alpha$ given by 
$\alpha(x_1,\dots,x_d) := (\alpha_1(x_1),\dots,\alpha_d(x_d))$ is an automorphism of 
$\Gamma$ and extends
$\alpha$. 

To prove that $\Gamma$ is Ramsey,
we show the statement for $d=2$; the general case then follows by induction on $d$. 
Let $\bA,\bB$ be substructures of 
$\Gamma = \Gamma_1 \boxtimes \Gamma_2$ and $r \in \mN$ be arbitrary. 
We will show that $\Gamma \to (\bB)^\bA_r$, 
so let 
$\chi \colon {\bB \choose \bA} \to [r]$ 
be arbitrary. If ${\bB \choose \bA}$ is empty,
then the statement is trivial, so in the following we assume that $\bA$ embeds into $\bB$.
For $i \in \{1,2\}$, let $\bA_i$ be the structure induced in $\Gamma_i$ by $\{a_i : (a_1,a_2) \in A\}$, and define $\bB_i$ analogously with $B$ instead of $A$.
Since $\Gamma_2$ is Ramsey there exists
a finite substructure $\bC_2$ of $\Gamma_2$
such that $\bC_2 \to (\bB_2)^{\bA_2}_r$.
Define $s := |{\bC_2 \choose \bA_2}|$.
Since $\Gamma_1$ is Ramsey there exists
a finite substructure $\bC_1$ of $\Gamma_1$
such that $\bC_1 \to (\bB_1)^{\bA_1}_{r^s}$. 
We identify the elements of $[r^s]$ with 
functions from ${\bC_2 \choose \bA_2}$ to $[r]$. 
Define $\chi_1 \colon {\bC_1 \choose \bA_1} \to [r^s]$ as follows. 
Let $e_1 \in {\bC_1 \choose \bA_1}$, 
let $e_2 \in {\bC_2 \choose \bA_2}$,
and let $e \in {\bB \choose \bA}$ be 
the embedding such that $e(a_1,a_2) = (e_1(a_1),e_2(a_2))$. Let $\xi \colon {\bC_2 \choose \bA_2} \to [r]$ be the function that maps
$e_2 \in {\bC_2 \choose \bA_2}$ to $\chi(e)$. 
Define $\chi_1(e_1) = \xi$. 
Then there exists an $f_1 \in {\bC_1 \choose \bB_1}$ such that 
$\chi_1(f_1 \circ {\bB_1 \choose \bA_1}) = \{\chi_2\}$
for some $\chi_2 \in {\bC_2 \choose \bA_2} \to [r]$.
As $\bC_2 \to (\bB_2)^{\bA_2}_r$, 
there exists an $f_2 \in {\bC_2 \choose \bB_2}$ 
such that 
$|\chi_2(f_2 \circ {\bB_2 \choose \bA_2})| = 1$. 
Let $f \in {\bC_1 \boxtimes \bC_2 \choose \bB}$ 
be given by $b \mapsto (f_1(b),f_2(b))$. 

We claim that 
$|\chi(f \circ {\bB \choose \bA})| = 1$. 
Arbitrarily choose $e,e' \in {\bB \choose \bA}$.
Then there are $e_i,e_i' \colon {\bB_i \choose \bA_i}$ for $i \in \{1,2\}$ such that 
$e(A) \subseteq (e_1(A_1),e_2(A_2))$
and $e'(A) \subseteq (e'_1(A_1),e'_2(A_2))$. 
Then $\chi_1(f_1 \circ e_1) = \chi_1(f_1 \circ e_1') = \chi_2$,
and $\chi_2(f_2 \circ e_1) = \chi_2(f_2 \circ e_2')$. 
Then $\chi(e) = \chi_2(f_2 \circ e_2) = \chi_2(f_2 \circ e_2') = \chi(e')$,
which is what we had to show. 
\end{proof}

The special case of Proposition~\ref{prop:product-ramsey} where
$\Gamma_1 = \cdots = \Gamma_d = ({\mathbb Q};<)$ can be found in~\cite{GrahamRothschildSpencer} (page 97).
The general case can also be shown inductively, see e.g.~\cite{BPT-decidability-of-definability}. 
One may also derive it using the results in Kechris-Pestov-Todorcevic~\cite{Topo-Dynamics}, since the direct product of extremely amenable groups is extremely amenable (also see~\cite{Bod-New-Ramsey-classes}). 

\subsection{Interpretations}\label{sect:interpret}
The concept of \emph{first-order interpretations} is a powerful tool to construct
new structures.
A simple example of an interpretation is
the line graph of a graph $G$,
which has a first-order interpretation over $G$.
%; a typical example 
%is the structure $({\mathbb Q};+,*)$ which has a first-order interpretation in $({\mathbb Z};+,*)$. 
By passing to the age of the constructed structure,
they are also a great tool to define new \emph{classes} of structures. 

\begin{Def}
A relational $\sigma$-structure $\bB$ has a \emph{(first-order) interpretation $I$} in a $\tau$-structure $\bA$ if there exists a natural number $d$, called the \emph{dimension} of $I$, and
\begin{itemize}
\item a $\tau$-formula $\delta_I(x_1, \dots, x_d)$ -- called the \emph{domain formula},
\item for each atomic $\sigma$-formula $\phi(y_1,\dots,y_k)$ a $\tau$-formula $$\phi_I(y_{1,1},\dots,y_{1,d}, y_{2,1}, \dots, y_{2,d},\dots,y_{k,1},\dots, y_{k,d})$$ -- the \emph{defining formulas};
\item a surjective map $h$ from
$\{\bar a : \bA \models \delta_I(\bar a)\}$ 
to $B$ -- called the \emph{coordinate map},
\end{itemize}
such that for all atomic $\sigma$-formulas $\phi$ and all elements $a_{1,1},\dots,a_{k,d}$ with $\bA \models \delta_I(a_{i,1},\dots,a_{i,d})$ for all $i \leq k$ 
\begin{align*}
& \bB \models \phi(h(a_{1,1},\dots,a_{1,d}), \dots, h(a_{k,1},\dots,a_{k,d})) \;  \\
 \Leftrightarrow \quad & 
\bA \models \phi_I(a_{1,1}, \dots,a_{k,d}) \; .
\end{align*}
\end{Def}

We give 
illustrating examples. 

\begin{example}
When $(V;E)$ is an undirected graph, then the \emph{line graph} of $(V;E)$ is the undirected graph $(E;F)$
where $F := \big \{\{u,v\} : |u \cap v| = 1 \big \}$. 
Undirected graphs can
be seen as structures where the signature contains a single binary
relation denoting a symmetric irreflexive relation. 
Then the line graph of $(V;E)$ 
has the following 2-dimensional interpretation 
$I$ over $(V;E)$: the domain formula $\delta_I(x_1,x_2)$ is $E(x_1,x_2)$, the defining formula for the atomic formula $y_1=y_2$ is $$(y_{1,1}=y_{2,1} \wedge y_{1,2} = y_{2,2}) \vee (y_{1,1}=y_{2,2} \wedge y_{1,2} = y_{2,1}) \, ,$$  and
the defining formula for the atomic formula
$F(y_1,y_2)$ is 
\begin{align*}
& \big((y_{1,1} \neq y_{2,1} \wedge y_{1,1} \neq y_{2,2}) \vee (y_{1,2} \neq y_{2,2} \wedge y_{1,2} \neq y_{2,2})\big) \\
\wedge & \; (y_{1,1}=y_{2,1} \vee y_{1,1} = y_{2,2} \vee y_{1,2} = y_{2,1} \vee y_{1,2} = y_{2,2}) \; .
\end{align*}
The coordinate map is the identity. 
\end{example}

\begin{example}
A poset $(P;\leq)$ has \emph{poset dimension at most $k$} if there are $k$ linear extensions
$\leq_1,\dots,\leq_k$ of $\leq$ such that 
$x \leq y$ if and only if $x \leq_i y$ for all $i \in [k]$. The class of all finite posets of poset dimension at most $k$ is the age of $({\mathbb Q};\leq)^k$, which clearly has a $k$-dimensional interpretation in $({\mathbb Q};<)$. 
\end{example}

\begin{lemma}[Theorem 7.3.8 in~\cite{HodgesLong}] 
\label{lem:interpret}
Let $\bA$ be an $\omega$-categorical structure. Then every
structure $\bB$ that is first-order interpretable in $\bA$ is countably infinite $\omega$-categorical or finite.
\end{lemma}

Note that in particular all reducts (defined in the introduction) of an $\omega$-categorical structure $\Gamma$ have an interpretation in $\Gamma$ and are thus again $\omega$-categorical. On the other hand, 
being homogeneous with finite relational
signature is not inherited by the interpreted structures. An example of a structure which is not interdefinable with a homogeneous structure in a finite relational signature, but which has a first-order interpretation over $({\mathbb N};=)$,
has been found by Cherlin and Lachlan~\cite{CherlinLachlan}. 

% Unclear: if $\Gamma$ is homogeneous
% in finite relational signature, can we have
% Delta also be homogeneous in finite
% relational signature?

\begin{prop}\label{prop:ramsey-interpret}
Suppose that $\Gamma$ is $\omega$-categorical Ramsey. 
Then every structure 
with a first-order interpretation in $\Gamma$
has an $\omega$-categorical Ramsey expansion $\Delta$. Furthermore, if $\Gamma$ is homogeneous with a finite relational signature,
then we can choose $\Delta$ to be homogeneous in a finite relational signature, too. 
\end{prop}

\begin{cor}
Conjecture~\ref{conj:ramsey-expansion} is true
for countable stable\footnote{For the definition of stability
we refer to any text book in model theory.} homogeneous structures with finite
relational signature. 
%All $\omega$-stable structures that
%are homogeneous in a finite relational signature can be expanded by finitely many relations to a homogeneous 
\end{cor}
\begin{proof}
Lachlan~\cite{LachlanIndiscernible} proved that every stable homogeneous
structure with a finite relational signature
has a first-order interpretation over $({\mathbb Q};<)$. The statement follows from the fact that
$({\mathbb Q};<)$ is Ramsey, and Proposition~\ref{prop:ramsey-interpret}. 
\end{proof}

\subsection{Adding constants}
\label{sect:constants}
% Why do we need constant symbols,
% rather than unary relations? It matters for
% the homogeneity business: consider
% for instance the random graph. When
% we add a singleton relation {r}, the resulting
% structure is no longer homogeneous,
% since we can map a point p such that E(r,p)
% to a point q such that N(q,p). 
% But we could avoid that by adding
% for each relation $R$ all relations
% defined by 
% $\exists p. (P(p) \wedge R(p,\bar x)). 
Let $\Gamma$ be homogeneous. It is clear that the expansion $(\Gamma,d_1,\dots,d_n)$ by finitely many constants $d_1,\dots,d_n$ is again homogeneous. Similarly, if
$\Gamma$ is $\omega$-categorical, then $(\Gamma,d_1,\dots,d_n)$ is $\omega$-categorical, as
a consequence of Theorem~\ref{thm:ryll}.
We will show here that if $\Gamma$ is Ramsey, then $(\Gamma,d_1,\dots,d_n)$ remains Ramsey. 
The original proof~\cite{BPT-decidability-of-definability}
went via a more general 
fact from topological dynamics (open subgroups of extremely amenable groups are extremely amenable). We
give an elementary proof here, due to Miodrag Sokic.

\begin{theorem}\label{thm:cramsey}
Let $\Gamma$ be homogeneous and Ramsey. 
Let $d_1,\dots,d_n$ be elements of $\Gamma$.
Then $(\Gamma,d_1,\dots,d_n)$ is also Ramsey.
\end{theorem}
\begin{proof}
Let $\tau$ be the signature,
and $D$ the domain of $\Gamma$. 
%and let $c_1,\dots,c_n$ be new constant symbols. 
%Let $\Gamma^*$ be an expansion of $\Gamma$ with
%signature $\tau \cup \{c_1,\dots,c_n\}$. 
We write
$d$ for $(d_1,\dots,d_n)$. 
Let $\bA^*,\bB^*$ be two finite
substructures of $\Gamma^*$, let $r \in \mN$, 
 and let $\chi^* \colon {\Gamma^* \choose \bA^*} \to [r]$ be arbitrary. We have to show that 
 there exists an $f \in {\Gamma^* \choose \bB^*}$
 such that $|\chi(f \circ {\bB^* \choose \bA^*})| = 1$.
We write $\bA$ and $\bB$ for the $\tau$-reducts of $\bA^*$ and $\bB^*$, respectively. 

Define $\chi \colon {\Gamma \choose \bA} \to [r]$
as follows. 
First, we fix for each tuple $a \in D^n$ 
that lies in the same orbit as $d$ 
in $\Aut(\Gamma)$ 
an automorphism $\alpha_a$ of $\Gamma$ such that
$\alpha_a(a) = d$. 
Let $e \in {\Gamma \choose \bA}$,
and $a := e(d)$. 
By the homogeneity of $\Gamma$, 
the tuples $a$ 
and $d$ lie in the same orbit 
of $\Aut(\Gamma)$. 
Note that $\alpha_{a} \circ e$ fixes $d$ and 
is an embedding
of $\bA^*$ into $\Gamma^*$. 
%We write $A^*_0$ for the copy of $A^*$ induced by 
%$\alpha_{\bar a}(e(A))$ in $\Gamma_c$. 
%Define $\chi(e(A)) := \chi^*(A_0^*)$.
Define $\chi(e) := \chi^*(\alpha_a \circ e)$. 
%Note that this is well-defined, since
%$e$ is the only embedding of $A$ into the
%substructure of $\Gamma$ induced by $e(A)$;
%otherwise, $A$ would not be rigid in contradiction to
%Corollary~\ref{cor:rigid} and $\Gamma$ being Ramsey. Brauchen wir nicht im setting
% wo wir embeddings faerben

Since $\Gamma$ is Ramsey, there is an $f \in {\Gamma \choose \bB}$
such that $\chi(f \circ {\bB \choose \bA}) = \{c\}$ for some $c \in [r]$. 
Let $b$ be $f(d)$.
By the homogeneity of $\Gamma$, 
the tuples $b$ 
and $d$ lie in the same orbit 
of $\Aut(\Gamma)$. 
Observe that 
$f' := \alpha_b \circ f$ fixes $d$ and is an embedding of $\bB^*$ into $\Gamma^*$.
%; write
%$B^*_0$ for the structure induced by
%$\beta(f(B))$ in $\Gamma_c$. 

%We claim that $\chi^*$ 
%is constant on $B^*_0 \choose A^*$. 
We claim that 
$|\chi^*(f' \circ {\bB^* \choose \bA^*})| = 1$.
To prove this, let $g \in {\bB^* \choose \bA^*}$ be arbitrary. 
%$g$ be an embedding of $A^*$ into $B^*_0$,
%and let $A^*_0$ be the copy of $A^*$ induced
%Then $g' := \alpha_{b}^{-1} \circ f' \circ e$ is an embedding of 
%$A$ into $\Gamma[f(B)]$, 
Since $g$ is in particular from ${\bB \choose \bA}$
we have $\chi(f \circ g) = c$.
%Since $\alpha_{b}(g'(A))$ induces $A_0^*$
%in $\Gamma_c$, 
By the definition of $\chi$ we have that
$\chi(f \circ g) = \chi^*(\alpha_b \circ f \circ g) = \chi^*(f' \circ g)$. 
Hence, $\chi^*(f' \circ g) = c$, which proves the claim.
\end{proof}

In this article, 
we work mostly with relational signatures. 
It is therefore important to note
that the relational structure $(\Gamma,\{d_1\},\dots,\{d_n\})$ is in general \emph{not} homogeneous even if $\Gamma$ is. Consider
for example the \Fresse\ limit $\Gamma = (\mV;E)$
of the class of all finite graphs, and an arbitrary $d_1 \in \mV$. Let $p \in \mathbb V \setminus \{d_1\}$ be such that $E(p,d_1)$ and $q \in \mathbb V  \setminus \{d_1\}$ be such that $\neg E(d_1,q)$. Then the mapping that sends $p$ to $q$ is an isomorphism between (one-element) substructures of $(\Gamma,\{d_1\})$ which cannot be extended to an automorphism of $(\Gamma,\{d_1\})$. (The difference to $(\Gamma,d_1)$ is that 
all substructures of $(\Gamma,d_1)$ must contain $d_1$.) 
Note, however, 
that $(\Gamma,d_1,\dots,d_n)$ and
$(\Gamma,\{d_1\},\dots,\{d_n\})$ have the
same automorphism group. 

%Recall that by Proposition~\ref{prop:ramsey-equiv}, the Ramsey properties of the 
%homogeneous expansion
%of an $\omega$-categorical structure 
%by all first-order definable relations only depends on the automorphism group of the structure. 

The solution to stating the result about expansions of homogeneous structures 
with constants in the relational
setting is linked to the following definition.  

\begin{Def}
Let $\Gamma$ be a relational
structure with signature $\tau$, and $d_1,\dots,d_n$ elements of $\Gamma$. 
Then $\Gamma_{d_1,\dots,d_n}$ 
denotes the expansion of $\Gamma$
which contains for every $R \in (\tau \cup \{=\})$
of arity $k \geq 2$, every $i \in [k]$ and $j \in [n]$, the $(k-1)$-ary relation $\{(x_1,\dots,x_{i-1},x_{i+1},\dots,x_k) : (x_1,\dots,x_k) \in R \text{ and } x_i = d_j\}$. 
\end{Def}

Note that if the signature of $\Gamma$ is finite,
then the signature of $\Gamma_{d_1,\dots,d_n}$ is also finite, and the maximal arity is unaltered. 
Also note that $\Gamma_{d_1,\dots,d_n}$ has in particular the unary relations $\{d_1\},\dots,\{d_n\}$.
%\item $\Gamma_{c_1,\dots,c_n}$ and 
%$(\Gamma,\{c_1\},\dots,\{c_n\})$ have the same automorphism group. Hence, when
%$\Gamma$ is $\omega$-categorical,
%%by Proposition~\ref{prop:ramsey-equiv}
%the expansion of 
%$\Gamma_{c_1,\dots,c_n}$ by all first-order 
%definable relations is Ramsey if and only if
%the expansion of 
%$(\Gamma,\{c_1\},\dots,\{c_n\})$ by all first-order definable relations is Ramsey. 
%\end{itemize}

\begin{lemma}
Let $\Gamma$ be a homogeneous relational
structure, and $d_1,\dots,d_n$ elements of $\Gamma$. Then $\Gamma_{d_1,\dots,d_n}$ is homogeneous. 
\end{lemma}
\begin{proof}
Let $a$ be an isomorphism between
two finite substructures $A_1, A_2$ 
of $\Gamma_{d_1,\dots,d_n}$. 
Since $\Gamma_{d_1,\dots,d_n}$ contains
for all $i \leq n$ the relation $\{d_i\}$ which is preserved by $a$, it follows that
if $A_1$ or $A_2$ contains $c_i$, then
both $A_1$ and $A_2$ must contain 
$d_i$, and $a(d_i)=d_i$. If $d_i$ is contained in neither $A_1$ nor $A_2$, then $a$ can be extended to a partial isomorphism $a'$ of
$\Gamma_{d_1,\dots,d_n}$ with domain $A_1 \cup \{d_i\}$ by
setting $a(d_i)=d_i$: this follows directly
from the definition of the signature of $\Gamma_{d_1,\dots,d_n}$. By the homogeneity of $\Gamma$, the map $a'$ can be extended to an automorphism of $\Gamma$. This
automorphism fixes $d_1,\dots,d_n$ pointwise,
and hence is an automorphism of $\Gamma_{d_1,\dots,d_n}$. 
\end{proof}

\begin{cor}
Let $\Gamma$ be homogeneous, $\omega$-categorical, and Ramsey,
and let $d_1,\dots,d_n$ be elements of $\Gamma$.
Then $\Gamma_{d_1,\dots,d_n}$ is also Ramsey.
\end{cor}
\begin{proof}
The statement follows from
Theorem~\ref{thm:cramsey}
from the observation that $\Gamma_{d_1,\dots,d_n}$
and $(\Gamma,d_1,\dots,d_n)$ 
have the same automorphism group, and
that whether an $\omega$-categorical structure
has the Ramsey property only depends
on its automorphism group (Proposition~\ref{prop:group}). 
\end{proof}

\subsection{Passing to the model companion}
\label{sect:mc}
A structure $\Gamma$ 
is called \emph{model-complete}
if all embeddings between models of the first-order theory of $\Gamma$ preserve all
first-order formulas. It is well-known that this is equivalent to every first-order formula being equivalent to an existential formula over $\Gamma$ (see e.g.~\cite{Hodges}). 
It is also known (see Theorem 3.6.7 in~\cite{Bodirsky-HDR-v4})  
that an $\omega$-categorical
structure $\Gamma$ 
is model-complete if and only if
for every finite tuple $t$ of elements of 
$\Gamma$ and for every self-embedding
$e$ of $\Gamma$ into $\Gamma$ there
exists an automorphism $\alpha$ of $\Gamma$
such that $e(t) = \alpha(t)$. 

%It follows from Proposition~\ref{prop:qe}
%that every homogeneous $\omega$-categorical structure is model-complete. 

A \emph{model companion} of $\Gamma$ is a 
model-complete structure $\Delta$ with the same age as $\Gamma$. 
If $\Gamma$ has a model companion, then
the model companion is unique up to isomorphism~\cite{Hodges}. Every $\omega$-categorical structure has a model 
companion, and the model companion is again
$\omega$-categorical~\cite{Saracino}. 

\begin{example}
We write $\mQ^+_0$ for $\{q \in \mQ : q \geq 0\}$. The structure $\Gamma := (\mQ^+_0;<)$ is $\omega$-categorical, but not model-complete: 
for instance the map $x \mapsto x+1$ is an embedding of
$\Gamma$ into $\Gamma$ which does not
preserve the unary relation $\{0\}$ with the 
first-order definition $\forall y (y \geq x)$ over $\Gamma$. The model companion of $\Gamma$ is $(\mQ;<)$. 
\end{example}

%An $\omega$-categorical structure $\Gamma$ has a model-complete theory if and only if all self-embeddings $e$
%of $\Gamma$ are locally generated by the automorphisms of
%$\Gamma$, that is, for every finite tuple $t$ of elements from $\Gamma$
%there exists an automorphism $\alpha$ of $\Gamma$ such that $e(t) = \alpha(t)$ (see e.g.~Theorem 3.6.11 in~\cite{Bodirsky-HDR}).
%In this case, we say that $\Gamma$ is model-complete. Note that by the above, when $\Gamma$ is homogeneous, 
%then $\Gamma$
%is model-complete. 
In this subsection we prove the following.

\begin{theorem}\label{thm:mc-Ramsey}
Let $\Gamma$ be $\omega$-categorical and Ramsey, 
and let $\Delta$ be 
the model companion of $\Gamma$. 
Then $\Delta$ is also Ramsey. 
\end{theorem} 
\begin{proof}
Let $e$ be an embedding of $\Gamma$ into $\Delta$,
and let $i$ be an embedding of $\Delta$ into $\Gamma$;
such embeddings exist by $\omega$-categoricity of $\Delta$ and $\Gamma$, see Section 3.6.2 in~\cite{Bodirsky-HDR-v4}. 
We will work with the equivalent 
characterisation of the Ramsey property given
in item 2 of Proposition~\ref{prop:group}.

Let $S$ and $M$ be finite subsets of the domain  $D$ of $\Delta$ and $r \in \mN$, and let
$\chi \colon {D \choose S} \to [r]$ be arbitrary. 
Let $D'$ be the domain of $\Gamma$. 
%Let $L$ be the set which exists 
We define a map $\chi' \colon {D' \choose i(S)} \rightarrow [r]$ as follows. For $q' \in {D' \choose i(S)}$, note that $e \circ q' \circ i \in {D \choose S}$. We define $\chi'(q') := \chi(e \circ q' \circ i)$.

Since $\Gamma$ is Ramsey, 
there exists
an $f' \in {D' \choose i(M)}$ and $c \in [r]$
such that for all $g' \in {i(M) \choose i(S)}$ 
we have $\chi'(f' \circ g') = c$.
Let $\alpha' \in \Aut(\Gamma)$ be an extension of $f'$. Note that $e \circ \alpha' \circ i$ is an embedding of $\Delta$ into $\Delta$, and since
$\Delta$ is model-complete there exists an $\alpha \in \Aut(\Delta)$ that extends the restriction $f$ of $e \circ \alpha' \circ i$ to $M$.

Let $g \in {M \choose S}$ be arbitrary.
We claim that $\chi(f \circ g) = c$. 
Since $e \circ i$ is an embedding of $\Delta$
into $\Delta$ and $\Delta$ is model-complete,
there exists an automorphism $\beta$ of $\Delta$ such that $\beta(e(i(x))) = x$ for all $x \in S$. Note that $g' := i \circ g \circ \beta \circ e \in {i(M) \choose i(S)}$, and hence $\chi'(f' \circ g') = c$. 
Also note that by the definition of $\chi'$ we have
\begin{align*}
\chi'(f' \circ g') = \chi(e \circ f' \circ g' \circ i) = \chi(e \circ f' \circ i \circ g \circ \beta \circ e \circ i) = \chi(f \circ g) \; .
\end{align*}
Hence, $\chi(f \circ g) = c$, and
$|\chi(f \circ {M \choose S})| \leq 1$, and thus $\Delta$ is Ramsey. 
\end{proof}

% An example of an application of Theorem~\ref{thm:mc-Ramsey} can be found at the end of Section~\ref{sect:inductive}.
% NOT IMPLEMENTED

\subsection{Passing to the model-complete core}
\label{sect:mc-core}
Cores play an important role in finite combinatorics. 
The concept of model-complete cores can be seen
as an existential-positive analog of model-companions, where embeddings 
are replaced by homomorphisms and self-embeddings are replaced
by endomorphisms. 
We state here results that are analogous to the results
for model companions that we have seen in the previous section. 
%We mention that the results for 
%model-complete cores of Ramsey structures 
%have been used in the p
%will be used
%in Section~\ref{sect:superimposing}. 

\begin{Def}
Let $\bA$ and $\bB$ be two structures with
domain $A$ and $B$, respectively, and
the same relational signature $\tau$. Then a 
 \emph{homomorphism} from $\bA$ to $\bB$
 is a function $f \colon A \to B$ such that for all
 $(a_1,\dots,a_n) \in R^\bA$ we have $(f(a_1),\dots,f(a_n)) \in R^\bB$.  
An \emph{endomorphism} of a structure $\Gamma$ 
is a homomorphism from $\Gamma$ to $\Gamma$. A structure $\Gamma$ is called 
a \emph{core} if every endomorphism of $\Gamma$ is
an embedding. 
\end{Def}

An $\omega$-categorical
structure $\Gamma$ 
is a model-complete core if and only if
for every finite tuple $t$ of elements of 
$\Gamma$ and for every endomorphism
$e$ of $\Gamma$ there
exists an automorphism $\alpha$ of $\Gamma$
such that $e(t) = \alpha(t)$ (Theorem 3.6.11 in~\cite{Bodirsky-HDR-v4}).  
%Hence, an $\omega$-categorical structure
%is a model-complete core if and only if the automorphisms
%of $\Gamma$ are dense in the endomorphisms of $\Gamma$.
The following has been shown in~\cite{Cores-Journal} (%for a more conceptual proof, 
also see~\cite{BodHilsMartin-Journal}).
Two structures $\Gamma$ and $\Delta$ are \emph{homomorphically equivalent} if there is a homomorphism from $\Gamma$ to $\Delta$ and a homomorphism from $\Delta$ to $\Gamma$.

\begin{theorem}\label{thm:cores}
Every $\omega$-categorical structure is homomorphically equivalent to
a model-complete core $\Delta$, which is unique up to isomorphism, and again countably infinite $\omega$-categorical or finite.
%For all $n$, orbits of $n$-tuples in $\omega$-categorical model-complete cores $\Gamma$ are primitive positive definable in $\Gamma$. %  DON'T NEED THIS: the
% paper is at most about ep properties, but never about pp 
% properties.
The expansion of $\Delta$ by all existential positive definable relations is homogeneous. 
\end{theorem}

The structure $\Delta$ in Theorem~\ref{thm:cores}
will be called \emph{the model-complete core of $\Gamma$}.

\begin{theorem}\label{thm:mc-core-Ramsey}
Let $\Gamma$ be $\omega$-categorical and Ramsey, 
and let $\Delta$ be the model-complete core of $\Gamma$. 
Then $\Delta$ is also Ramsey. 
\end{theorem} 
\begin{proof}
The proof is similar to the proof of Theorem~\ref{thm:mc-Ramsey}. 
\end{proof}
%An example for an application of Theorem~\ref{thm:mc-core-Ramsey} will be given in Section~\ref{sect:inductive}. 

\subsection{Superimposing signatures}
\label{sect:superimposing}
An amalgamation class $\mathcal C$ 
is called a \emph{strong amalgamation class} if,
informally,  
we can amalgamate structures $\bB_1,\bB_2 \in \mathcal C$ over $\bA \in \mathcal C$ in such a way that no points of $\bB_1$ and $\bB_2$  other than the elements of $\bA_1$ will be
identified in the amalgam. Formally, 
we require that
for all $\bA,\bB_1,\bB_2 \in \mathcal C$ and embeddings $e_i \colon \bA \to \bB_i$, $i \in \{1,2\}$,
there exists a structure $\bC \in \mathcal C$ and embeddings $f_i \colon \bB_i \to \bC$ such that
$e_1(f_1(x)) = e_2(f_2(x))$ for all $x \in A$,
and additionally $f_1(B_1) \cap f_2(B_2) = f_1(e_1(A)) = f_2(e_2(A))$. 
When an amalgamation class $\mathcal C$ even has
strong amalgamation, then this can be seen
from the automorphism group of the \Fresse\ limit of $\mathcal C$. 

\begin{Def}[\cite{Oligo}]
We say that a permutation group has \emph{no algebraicity} if for every finite tuple $(a_1,\dots,a_n)$ of the domain the set of all permutations of the group that fix each of $a_1,\dots,a_n$ fixes no other elements of the domain. 
\end{Def}

For automorphism groups of $\omega$-categorical structures $\Gamma$, 
having no algebraicity
coincides with the model-theoretic notion of
$\Gamma$ having no algebraicity (see, e.g.,~\cite{Hodges}). 

\begin{lemma}[see (2.15) in \cite{Oligo}]
Let $\C$ be an amalgamation class of relational structures and $\Gamma$ its \Fresse~limit. Then $\C$ has strong amalgamation if and only if $\Gamma$ has no algebraicity. 
\end{lemma}

For strong amalgamation classes there is a powerful construction to obtain new strong amalgamation classes from known ones. 

\begin{Def}
Let $\C_1$ and $\C_2$ be classes of finite structures with disjoint relational signatures $\tau_1$ and $\tau_2$, respectively. Then the
\emph{free superposition of $\C_1$ and $\C_2$}, denoted by $\C_1 * \C_2$, is the class of $(\tau_1 \cup \tau_2)$-structures $\bA$ such that the $\tau_i$-reduct of $\bA$ is in $\C_i$, for $i \in \{1,2\}$.
\end{Def}

The following lemma has a straightforward proof
by combining amalgamation in 
$\C_1$ with amalgamation in $\C_2$.

\begin{lemma}
If $\C_1$ and $\C_2$ are strong amalgamation classes, then $\C_1 * \C_2$
is also a strong amalgamation class. 
\end{lemma}
%\begin{proof}
%\end{proof}

When $\Gamma_1$ and $\Gamma_2$ are 
homogeneous structures with no algebraicity,
then $\Gamma_1 * \Gamma_2$
denotes the (up to isomorphism unique)
\Fresse\ limit of the free superposition of the
age of $\Gamma_1$ and the age of $\Gamma_2$. 

\begin{example}\label{expl:random-perm}
For $i \in \{1,2\}$, let $\tau_i = \{<_i\}$,
let $\C_i$ be the class of all finite $\tau_i$-structures where $<_i$ denotes a linear order,
and let $\Gamma_i$ be the \Fresse\ limit of $\C_i$. 
Then $\Gamma_1 * \Gamma_2$
is known as the \emph{random permutation} (see e.g.~\cite{BoettcherFoniok,LinmanPinsker,Sokic}). 
\end{example}

We have the following result about free superpositions.

\begin{theorem}[\cite{Bod-New-Ramsey-classes}]\label{thm:superposition}
Let $\Gamma_1$ and $\Gamma_2$ be
homogeneous $\omega$-categorical 
structures with no algebraicity
such that both $\Gamma_1$ and $\Gamma_2$
are Ramsey. Then $\Gamma_1 * \Gamma_2$ is Ramsey. 
\end{theorem}

We mention that the proof of Theorem~\ref{thm:superposition} from~\cite{Bod-New-Ramsey-classes} uses Theorem~\ref{thm:mc-core-Ramsey}
 about model-complete cores. An alternative proof can be found in~\cite{Sok-Directed-graphs-and}. 

\begin{example}\label{expl:E-with-two-orders}
Recall from Example~\ref{expl:equivalences}
that the amalgamation class of all finite structures $(V;E,<)$
where $E$ denotes an equivalence relation and
$<$ denotes a linear order, is \emph{not} Ramsey.
In the light of Conjecture~\ref{conj:ramsey-expansion} for the \Fresse\ limit $\Gamma$ of this class, we therefore look for a homogeneous Ramsey expansion of $\Gamma$. Let $\mathcal C$ be the class of
all finite structures $(V;E,\prec)$ where
$E$ is an equivalence relation and $\prec$
is a linear order that is convex with respect to $E$.
We have mentioned before that $\mathcal C$
is Ramsey, and by Theorem~\ref{thm:superposition} the class $\mathcal C * \LO$ is Ramsey. Then the \Fresse\ limit of 
$\mathcal C * \LO$ is isomorphic to a homogeneous
Ramsey expansion of $\Gamma$. 
\end{example}

\begin{example}\label{ex:s2}
The directed graph $S(2)$ is one of the 
homogeneous directed graphs that figures in the classification of all homogeneous directed graphs of Cherlin~\cite{Cherlin}. In fact, it is a homogeneous
tournament and therefore already appeared in the
classification of homogeneous tournaments of Lachlan~\cite{Lachlan}. 
It has many equivalent definitions, one of them
being the following: the vertices of $S(2)$
are a countable dense set of points on the unit circle without antipodal points. We add an edge from $x$ to $y$ if and only if the line from $x$ to $y$ has the origin on the left; that is, $x$, $y$, and $(0,0)$ lie
in clockwise order in the plane. 

We will show that $S(2)$ has a Ramsey expansion which is homogeneous and has a finite relational signature. This can be derived from general
principles and Ramsey's theorem as follows. 
In the following, $(\mQ;<_1)$ and $(\mQ;<_2)$
both denote the order of the rationals, but
have disjoint signature. 
Let $\Gamma$ be the disjoint union 
$({\mathbb Q};<_1) \uplus_P ({\mathbb Q};<_1)$, 
which has the Ramsey property 
by Theorem~\ref{thm:ramsey} and Lemma~\ref{lem:disj-union} (Example~\ref{expl:lo}).  
%It is easy to see that the age of $\Gamma$ has strong amalgamation. 
%Define $x \prec_1 y$ by $x < y \vee (P(x) \wedge \neg P(y))$. 
Then the free superposition $\Delta$
of $({\mathbb Q};<_2)$ with $\Gamma$ 
is Ramsey by Theorem~\ref{thm:superposition},
and homogeneous with finite relational signature. 
The structure $S(2)$ is a reduct of $\Delta$:
for elements $x,y \in S(2)$, 
we define $x \prec y$ if $$\big (x <_2 y \wedge (P(x) \Leftrightarrow P(y)\big) \vee \big(y <_2 x \wedge (P(x) \not\Leftrightarrow P(y)\big) \; .$$
\end{example}

Let $\Gamma_1$ and $\Gamma_2$ be two 
$\omega$-categorical Ramsey structures. 
Note that since there is a linear order with a first-order definition in $\Gamma_1$, and a first-order
definition of a linear order in $\Gamma_2$,
the structure
$\Gamma_1 * \Gamma_2$ must carry two independent
linear orders. 

To prove the Ramsey property for structures that
do not have a second independent linear order,
we have the following variant. 

\begin{theorem}[\cite{Bod-New-Ramsey-classes}]
Let ${\mathcal C}_1$ and ${\mathcal C}_2$
be classes of structures such that
${\mathcal C}_1$, ${\mathcal C}_2$
 and $\LO$ have pairwise disjoint signatures.
Also suppose that $\mathcal C_1$ and $\LO * {\mathcal C}_2$ are Ramsey classes with strong amalgamation and 
$\omega$-categorical \Fresse\ limits.
Then ${\mathcal C_1} * {\mathcal C}_2$
is also a Ramsey class. 
\end{theorem}

% 15/5/14: 
% There is a problem in Roedl-Nesetril: 
% in the partite Lemma, they produce triangles. 
% Same in Proemel Voigt?

% This is why I am using here an idea from
% Lionel: the idea is to work with "liftings":
% an $n$-partite structure B is a lifting 
% of an $n$-partite structures A if whenever
% (a1,...,an) is an edge of B, then (pi(a1),...,pi(an))
% is an edge of A. 

% Note that there is a difference to "rectified":
% we don't require "if and only if".
% Several advantages here: 
% 1) No triangles show
% up if they were not present before. 
% 2) In the partite amalgamation construction,
% we only need to apply the partite Lemma
% to structures that are lifts of A, WHEN
% WE START FROM A "Glueing-template"
% obtained from the old proof which
% creates triangles. BUT have to be careful:
% amalgamation on the right might destroy the
% assumption that we have a glueing-template,
% need some variant free amalgamation for this technique. 

% FUNNY: still need the old partite method to produce this template!

\section{The partite method} 
%of \Nesetril\ and R\"odl}
\label{sect:partite}
There are some homogeneous Ramsey structures 
where no proof of the Ramsey property from general
principles is known. One of the most powerful
methods to prove the Ramsey property is such situations 
is the \emph{partite method}. 
The first result that we see in this section 
is that for any finite relational signature $\tau$, 
the class of all finite ordered $\tau$-structures
is a Ramsey class. This is due to \Nesetril\ and R\"odl~\cite{NesetrilRoedlOrderedStructures} 
and independently to Abramson and Harrington~\cite{AbramsonHarrington}; in these original papers,
the statement is made for hypergraphs only,
but it holds for relational structures in general. 
We then apply the partite method to classes that are characterised by forbidding finite structures as induced substructures. 

\subsection{The class of all ordered structures}
\label{sect:all-structs}
We will prove the following theorem, due to \Nesetril\ and R\"odl, and, independently, Abramson and Harrington. 

\begin{theorem}[\cite{AbramsonHarrington,NesetrilRoedlOrderedStructures}]\label{thm:ordered-graph-ramsey}
For every relational signature $\tau$ the class
of all $\tau \cup \{\preceq\}$-structures, where $\preceq$ denotes a linear order, is a Ramsey class. 
%The class of all ordered digraphs $G = (V;E,<)$
%is a Ramsey class.
%For all 
%$A,B \in {\mathcal G}^<$ and 
%every $k > 1$ there exists a $C \in {\mathcal G}^<$
%such that $C \rightarrow (B)^A_k$. 
\end{theorem}

The construction to prove the Ramsey property is due to \Nesetril\ and R\"odl~\cite{NesetrilRoedlPartite}, 
with only minor modifications in the presentation.
%which we will describe after Definition~\ref{def:power}). 
%However, we provide more proof details here since we will later modify the construction.
% to exhibit some new Ramsey classes. 
%The proof technique of \Nesetril\ and R\"odl is called the \emph{partite method}, 
It relies on the concept of
\emph{$n$-partite structures}. 
%To facilitate the generalization
%to other classes of structures, our 
We formalize this slightly differently than 
\Nesetril\ and R\"odl in~\cite{NesetrilRoedlPartite}. %, in order to 

% Differences: 
%1 allow edges on the same level (check: really OK?)
% 2 Coding: use weak partial orders instead of multi-sorted domain and building everything into the notion of embeddings. 

\begin{Def}
Let $n \in \mathbb N$, and $\tau$ a relational signature. An \emph{$n$-partite structure} is a finite $(\tau \cup \{\preceq\})$-structure $(\bA,\preceq)$
%$(V;E,\preceq)$ where $E$
%is a binary relation and 
where $\preceq$ is a weak linear order (that is, a linear quasi-order)
such that the equivalence relation $\approx$ on $A$ defined by $x \approx y \Leftrightarrow (x \preceq y \wedge y \preceq x)$ has $n$ equivalence classes.
%\item $x \approx y \Rightarrow \neg E(x,y)$. 
%\end{itemize}
An $n$-partite structure is called a \emph{transversal} 
%\emph{transversal} 
if each equivalence
class of $\approx$ has size one. 
\end{Def}

Note that the elements of a finite $n$-partite structure $(\bA,\preceq)$  are partitioned into \emph{levels} $A_1,\dots,A_n$ which
are uniquely given by the property that for $u \in A_i$ 
and $v \in A_j$ we have $u \preceq v$ 
if and only if $i \leq j$.

\begin{lemma}[Partite Lemma]\label{lem:partite}
Let $\bA$ be an $n$-partite transversal, $\bB$ an arbitrary $n$-partite structure, and $r \in {\mathbb N}$. 
Then there exists an $n$-partite structure $\bC$ 
such that $\bC \rightarrow (\bB)^\bA_r$.
\end{lemma}

The idea of the proof of Lemma~\ref{lem:partite}
is to use the theorem of Hales-Jewett (see~\cite{GrahamRothschildSpencer}),
which we quickly recall here to fix some terminology. 

\begin{Def}
Let $m,d \in {\mathbb N}$.
A \emph{combinatorial line} is a 
set $L \subseteq [m]^d$ of the form 
$$\{(\alpha^1_1,\dots,\alpha^1_d),\dots,(\alpha^m_1,\dots,\alpha^m_d)\}$$ 
such that there exists a non-empty set
$P_L \subseteq [d]$ satisfying
\begin{itemize}
\item $\alpha^k_p = \alpha^l_p$ for all $k,l \in [m]$ and $p \in [d] \setminus P_L$, and
\item $\alpha^k_p = k$ for all $k \in [m]$ and $p \in P_L$. 
\end{itemize}
%$$ L = \big \{(\alpha_1,\dots,\alpha_d) \in [r]^d \; \big | \; \alpha_l =\alpha_l \text{ for } l \in [d] \setminus P_L, \text{ and } \alpha_l=\alpha_m \text{ for } l,m \in P_L \big \} \; .$$
\end{Def}
Note that for every $k \in [m]$ there exists exactly one
$\alpha = (\alpha_1,\dots,\alpha_d) \in L$ with $\alpha_p = k$ for all $p \in P_L$;
we write $L(k)$ for this $\alpha$. 

\begin{figure}
\begin{center}
  	\includegraphics[scale=0.3]{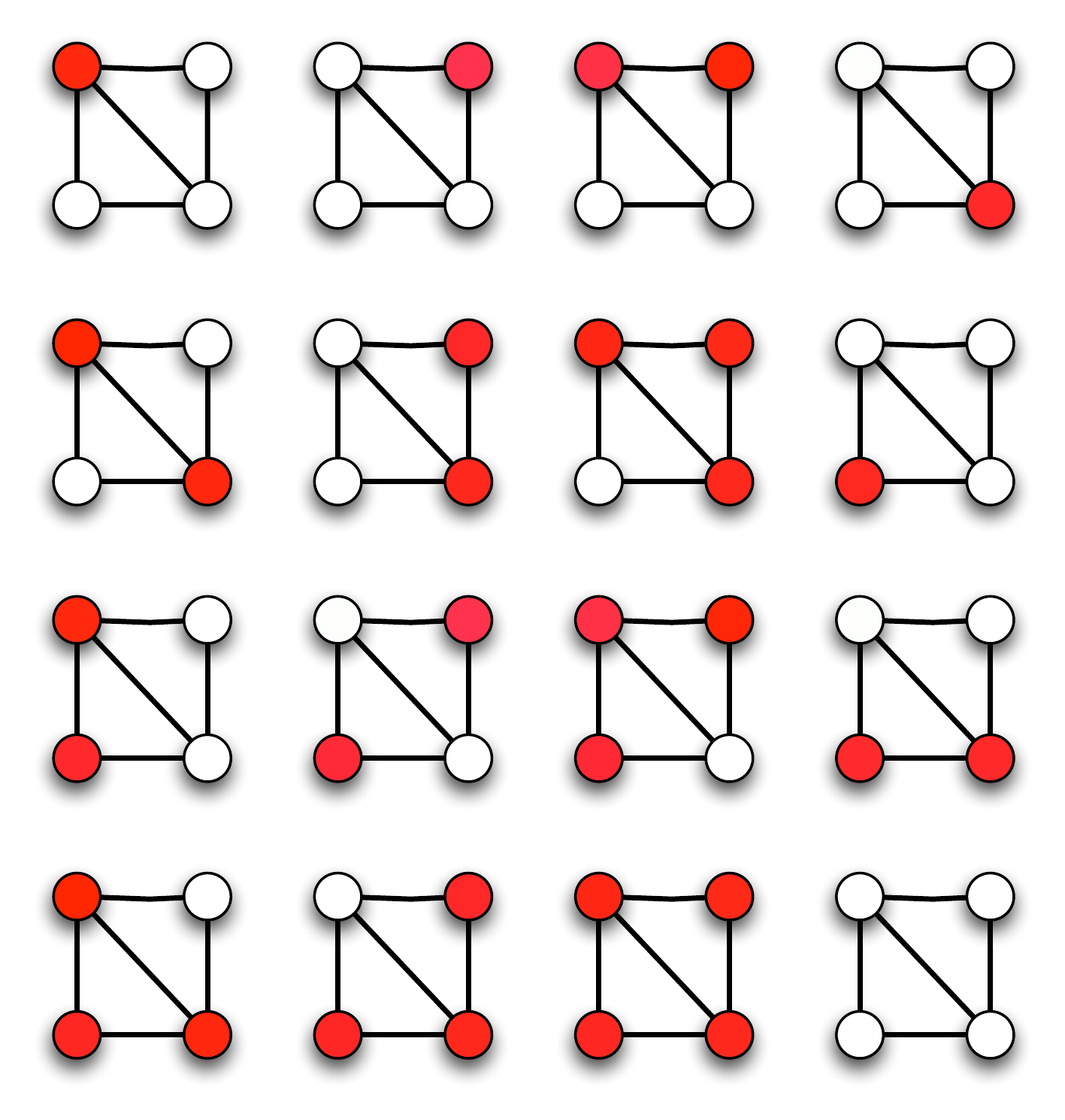} 
\caption{Illustration for $HJ(2,2)=2$.}
\end{center}
\label{fig:HJ}
\end{figure}

\begin{theorem}[Hales-Jewett; see~\cite{GrahamRothschildSpencer}]
\label{thm:HJ}
For any $m,r \in {\mathbb N}$  there exists $d \in {\mathbb N}$
such that for every function $\xi \colon [m]^d \to [r]$ 
there exists a combinatorial line $L$ such that
$\xi$ is constant on $L$. 
\end{theorem}
We write $HJ(m,r)$ 
for the smallest $d  \in {\mathbb N}$ 
that satisfies the condition in Theorem~\ref{thm:HJ}. 
See Figure~\ref{fig:HJ} for an illustration that shows
that $HJ(2,2) = 2$: if we colour the vertices of $[2]^2$ with two colours,
we always find a monochromatically coloured combinatorial line.

\begin{proof}[of Lemma~\ref{lem:partite}] 
We assume that every vertex of $\bB$
is contained in a copy of $\bA$ in $\bB$. This is without loss of generality:
if $\bB^*$ is the substructure of $\bB$ induced by the elements of the copies of $\bA$ in $\bB$,
and $\bC^*$ is such that $\bC^* \rightarrow (\bB^*)^{\bA^*}_r$, then we can construct $\bC$ such that $\bC \rightarrow (\bB)^\bA_r$ from $\bC^*$ by amalgamating at every copy of $\bB^*$ in $\bC^*$ a copy of $\bB$.  
So assume in the following that $\bB = \bB^*$. 

%TODO: more formally, please. 
% WHEN DOES THIS STEP WORK, IN GENERAL?
% It seems that we only need amalgamation here! 
% Don't have to work in the partite setting for this!
% Make it a lemma earlier on?

Let $g_1,\dots,g_m$ be an enumeration
of ${\bB \choose \bA}$. 
Let $d$ be $HJ(m,r)$ 
(according to Theorem~\ref{thm:HJ}). 
The idea of the construction in the proof of Lemma~\ref{lem:partite} is to construct $\bC$ in such a way
that for every element of $[m]^d$ 
there exists a copy of $\bA$ in $\bC$ such that 
monochromatically coloured lines in $[m]^d$ correspond
to monochromatic copies of $\bB$ in $\bC$. 
The direct product $\bB^d$ has many copies of $\bA$,
but in general does not have enough copies of $\bB$. 
The following ingenious construction, named after
the initials of its inventors, is a modification 
of the direct product that overcomes the  
mentioned problem by creating sufficiently many copies 
of $\bB$. 

\begin{Def}[The NR-power]\label{def:power}
Let $\bA$, $\bB$ be $n$-partite structures with signature $\tau$. 
Then the \emph{$d$-th NR-power of $\bB$ over $\bA$} is the $n$-partite structure $\bC$ defined as follows. 
Write $B_i$ for the $i$-th level of $B$, for $i \in [n]$. The domain of $\bC$ is
 $C_1 \cup \cdots \cup C_n$ where $C_i := (B_i)^d$. 
% For $u,v \in C$, define $u \preceq^\bC v$ iff $u_p \preceq v_p$ for all $p \in [d]$.
%Finally, let 
For 
$R \in \tau$ of arity $h$,
and $u^1,\dots,u^h \in C$, we define
$(u^1,\dots,u^h) \in R^\bC$ iff
%\begin{enumerate}
%\item $(u^1_p,\dots,u^h_p) \in R^\bB$ 
%for all $p \in [d]$, or
%\item 
\begin{itemize}
\item there is a non-empty set $P \subseteq [d]$ and $(w^1,\dots,w^h) \in R^\bB$ 
such that $u^s_q=w^s$
for $q \in P$ and $s \in [h]$, and 
\item for $q \in [d] \setminus P$, all of $u^1_q,\dots,u_q^h$ lie in the same copy of $\bA$ in $\bB$.
\end{itemize} 
\end{Def}

For an illustration of the NR-power, see Figure~\ref{fig:partite}. 

\begin{figure}
  	\includegraphics[scale=0.6]{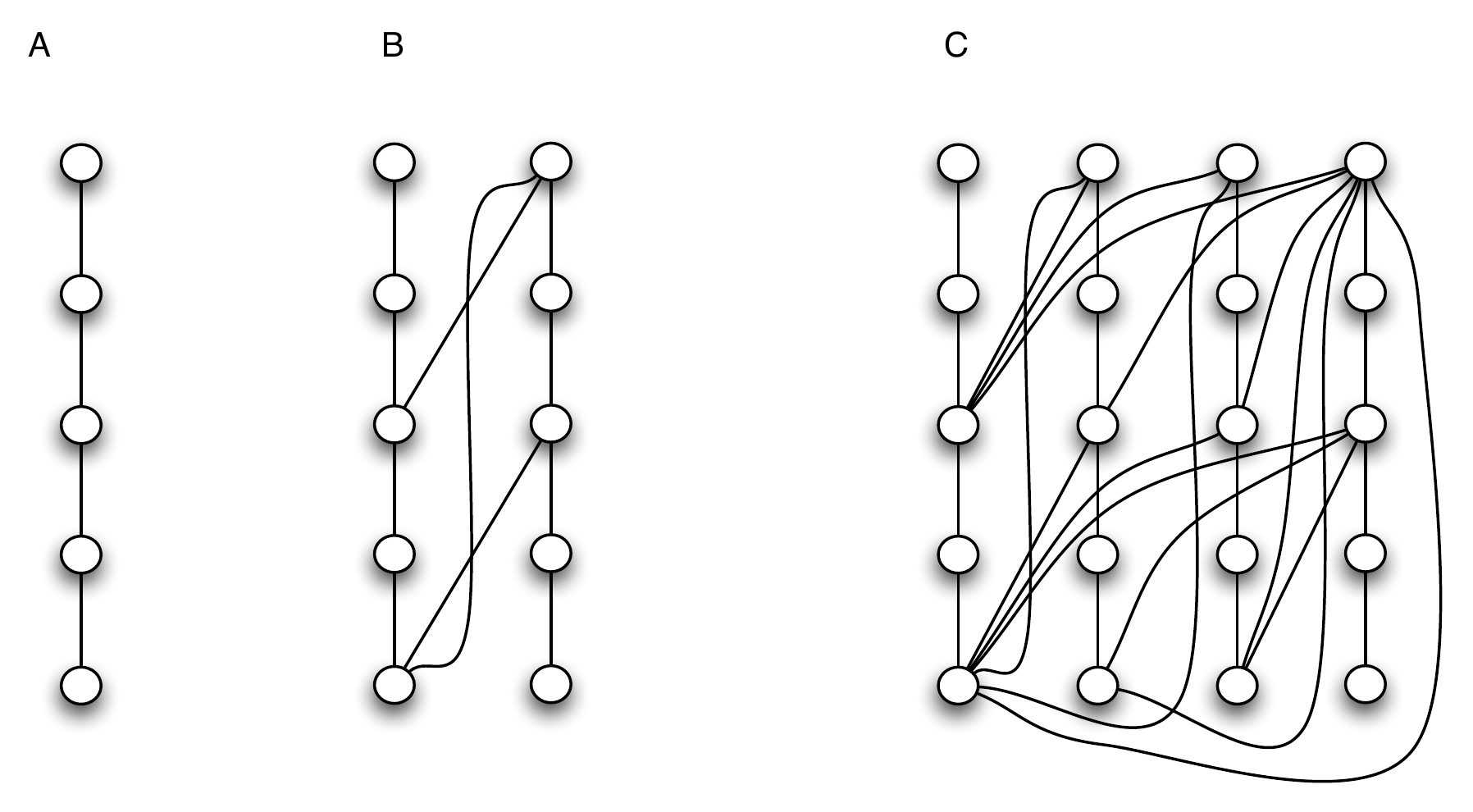} 
\caption{Illustration of two five-partite graphs $\bA$, $\bB$, where $\bA$ is transversal. On the right, we see the second NR-power of $\bB$ over $\bA$.}
\label{fig:partite}
\end{figure}

% Frage 1: Warum diese Guard am Schluss? 
% Antwort 8/5/14: 
% wir muessen ja verschiedene Sachen sicherstellen:
% 1) die linien muessen sich alle die gleichen Kopien von den
% A's teilen --> das spricht fuer "direct-product" Kopien von A's,
% und gegen Vorkommen von A's die "vertikal" oder "horizontal"
% im Box-product vorkommen wuerden.
% wenn einfach nur das strong product
% haetten wir nicht genuegend Kopien von 
% B's: Siehe das Beispiel mit den zwei Pfaden,
% wo es eine zusaetzliche Kante gibt vom ersten Knoten des ersten Pfades zum letzten
% Knoten des letzten Pfades. 
% 2) wir wollen wirklich Kopien von A's und B's kriegen
% Beispiel dazu: A = Kante, B = zwei disjunkte Kanten.
% Ohne guards kriegen wir unerwuenschte Querquanten 8/5/14. STIMMT 8/7/14

% Frage 2: 
% Warum denn immer die gleiche Kante? 1/5/14
% Antwort 8/7/14: Glaube es geht beides.
% Betrachte das Beispiel mit A = pfad, B = den zwei Pfaden 135 und 246, und den extra Kanten 
% 16 und 25. Wenn das nicht
% die gleiche Kanten sein muss, dann haetten 
% wir auf der Diagonalen vom HJ 
% folgendes: 
% Pfad 11,33,55, Pfad 22,44,66, und extra Kanten von 11 nach 66 und 22 nach 55: 
% GIBT ALSO KEIN PROBLEM

\medskip 

Let $\bC$ be
the $d$-th NR-power of $\bB$ over $\bA$.
To prove the partite lemma, it suffices to show that $\bC \to (\bB)^\bA_r$. Let $\chi \colon {\bC \choose \bA} \to [r]$ be arbitrary. We are going to define a function 
$\xi \colon [m]^d \to [r]$. 

%For $\alpha = (\alpha_1,\dots,\alpha_d) \in [m]^d$ and $a \in A$, we
%write $a^\alpha$ for the element $(g_{\alpha_1}(a),\dots,g_{\alpha_d}(a))$ of $C$.

\vspace{.2cm}
{\bf Claim 1.} 
For $\alpha = (\alpha_1,\dots,\alpha_d) \in [m]^d$, the map $g_\alpha \colon A \to C$ given by
$a \mapsto (g_{\alpha_1}(a),\dots,g_{\alpha_d}(a))$
%$a \mapsto a^\alpha$ 
is an embedding of $\bA$ into $\bC$. 

\vspace{.1cm}

%It is straightforward to see that $g_\alpha$ preserves $\preceq$: we have $a \preceq^\bA b$ if and only if $g_{\alpha_p}(a) \preceq^\bB g_{\alpha_p}(b)$ for all $p \in [d]$ 
%since $g_{\alpha_p}$ is an embedding. This in turn is equivalent, by definition of the partite product, to $g_\alpha(a) \preceq^\bC g_\alpha(b)$.
%Now 

\begin{proof}[of Claim 1.]
Suppose that 
$(a_1,\dots,a_h) \in R^{\bA}$. 
Then 
$$(g_{\alpha_p}(a_1),\dots,g_{\alpha_p}(a_h)) \in R^{\bB}$$ for all $p \in [d]$
since $g_{\alpha_p}$ preserves $R$.
By the definition of $R^\bC$ 
we have that $(g_\alpha(a_1),\dots,g_\alpha(a_h)) \in R^\bC$ (arbitrarily choose $i \in [d]$ and verify
Definition~\ref{def:power} for $P = \{i\}$).
Conversely, suppose that $(g_\alpha(a_1),\dots,g_\alpha(a_h)) \in R^\bC$.
Then there exists a non-empty set $P \subseteq [d]$ and $(w_1,\dots,w_h) \in R^\bB$ 
such that for all $q \in P$ and $s \in [h]$
we have $(g_\alpha(a_s))_q = g_{\alpha_q}(a_s) = w_s$. 
Since $g_{\alpha_q}$ is an embedding of $\bA$ into $\bB$, we obtain
in particular that $(a_1,\dots,a_h) \in R^\bA$, proving the claim. 
\end{proof}

Define $\xi(\alpha) := \chi(g_\alpha)$. 
By the theorem of Hales-Jewett (Theorem~\ref{thm:HJ}), 
there exists a combinatorial line $L \subseteq [m]^d$ and $c \in [r]$
such that $\xi(\alpha) = c$ for all $\alpha \in L$. 
We describe how $L$ gives rise to an
embedding $g_L$ of $\bB$ into $\bC$.
For $u \in B$, we write $\pi(u)$ for the unique element of $A$ that lies on the same level as $u$. Observe that $\pi(g_k(a)) = a$ for all $k \in [m]$ and for all $a \in A$ since $\bA$ is transversal.
Recall our assumption that 
every $u \in B$ appears in a copy of
$\bA$ in $\bB$, and hence 
there exists a $k \in [m]$ such that
$u \in g_k(A)$. 

\vspace{.2cm}
{\bf Claim 2.}
The map $g_L \colon B \to C$
given by $g_L(u) := g_{L(k)}(\pi(u))$, 
for some $k \in [m]$ such that $u \in g_k(A)$,
 is well-defined, and an embedding of $\bB$ into $\bC$. 
\vspace{.1cm}

\begin{proof}[of Claim 2.]
In order to show that the value 
of $g_{L}$ does not depend on the choice
of $k$, we have to show that if 
there are $k,l \in [m]$
such that $u \in B$ appears in both 
$g_k(A)$ and in $g_l(A)$, then 
$g_{L(k)}(\pi(u)) = g_{L(l)}(\pi(u))$,
that is, $g_{L(k)_p}(\pi(u))  = g_{L(l)_p}(\pi(u))$ for all $p \in [d]$. 
This is clear when $p \in [d] \setminus P_L$ since we then have $L(k)_p = L(l)_p$. 
So consider the case $p \in P_L$. 
Then
\begin{align*}
g_{L(k)_p}(\pi(u)) = g_k(\pi(u)) = u = g_l(\pi(u)) = g_{L(l)_p}(\pi(u)) 
\end{align*} 
where the equation $g_k(\pi(u)) = u = g_l(\pi(u))$ holds
since $A$ is transversal and $u \in g_k(A) \cap g_l(A)$. 

To show that $g_L$ is an embedding, 
let $R \in \tau$ be of arity $h$, 
and let $u_1,\dots,u_h \in B$ be arbitrary.
%First suppose that there exists a $k \in [m]$
%such that $u_1,\dots,u_h$ all lie in $g_k(A)$. 
%Note that $g_L(u_s) = g_{L(k)}(\pi(u_s))$ for 
%all $s \in [h]$. 
%Then
%\begin{align*}
%(u_1,\dots,u_h) \in R^\bB 
%\; \Leftrightarrow \; & (\pi(u_1),\dots,\pi(u_h)) \in R^{\bA} \\
%\Leftrightarrow \; & (g_{{L(k)}}(\pi(u_1)),\dots,g_{L(k)}(\pi(u_h))) \in R^{\bC} \\
%\Leftrightarrow \; & (g_{L}(u_1),\dots,g_{L}(u_h)) \in R^{\bC}
%\end{align*}
%where the middle equivalence is by Claim 1. 
%So suppose in the following that $u_1,\dots,u_h$ do not all lie in the same copy of $\bA$ in $\bB$. 
Let $s \in [h]$ and $k$ be such that $u_s \in g_k(A)$. Let $p \in P_L$ be arbitrary. Then
\begin{align}
(g_L(u_s))_p = (g_{L(k)}(\pi(u_s)))_p = g_{L(k)_p}(\pi(u_s)) = g_k(\pi(u_s)) = u_s  \; .
\label{eq:central-partite}
\end{align}
%%In particular, $((g_L(u_1))_p, \dots, (g_L(u_h))_p) = (u_1,\dots,u_s) \in R^\bB$. 
%In particular, 
%%$(g_L(u))_l = u$ and $(g_L(v))_l$ do not lie in the same copy of $A$ in $B$ for $l \in P_L$. 
%For $q \in [d] \setminus P_L$ we have
%$L(k)_q = L(l)_q$ for all $k,l \in [m]$, 
%by the properties of combinatorial lines. 
%Since $(g_L(u_s))_q = g_{L(k)_q}(\pi(u_s))$
%for some $k \in [m]$, 
%all of $(g_L(u_1))_q,\dots,(g_L(u_h))_q$
%lie in the same copy of $\bA$ in $\bB$.
Hence, if $(u_1,\dots,u_s) = ((g_L(u_1))_p, \dots, (g_L(u_h))_p) \in R^\bB$, then by the  definition of $R^\bC$ for $P := P_L$ and $w^s := u_s$ for all $s \in [h]$ we have that $(g_L(u_1),\dots,g_L(u_h)) \in R^\bC$, and $g_L$ preserves $R$.

Conversely, suppose that
$(g_L(u_1),\dots,g_L(u_h)) \in R^\bC$. 
Then there is a non-empty set $P \subseteq [d]$
and $(w^1,\dots,w^h) \in R^\bB$ such that
for $q \in P$ and $s \in [h]$ we have $g_L(u_s)_q = w^s$, and for $q \in [d] \setminus P$, 
all of $g_L(u_1)_q,\dots,g_L(u_h)_q$ lie in the same copy of $\bA$ in $\bB$. 
For $p \in P$ we have $w^s = (g_L(u_s))_p = u_s$, 
and thus $(u_1,\dots,u_h) \in R^\bB$. 
%Otherwise, $(g_L(u_1))_p=u_1,\dots,(g_L(u_s))_p=u_s$
%all lie in the same copy of $\bA$ in $\bB$ contrary to our assumptions. 
Applied to the case where $R$ is the equality
relation (for proving Ramsey results, 
we can assume without loss of generality that the signature contains a symbol for equality),
this also shows injectivity of $g_L$.
Hence, $g_L$ is an embedding, 
which concludes the proof of the claim. 
\end{proof}

%Now let $A' \in {B^L \choose A}$ be arbitrary. 
%By Claim 1, the structure induced by $A^\alpha \subseteq B^L$ induces in $C$ a structure isomorphic to $A$ for each $\alpha \in L$. Since $\xi(\alpha) = c_0$ for all $\alpha \in L$,
 %we also have $\chi(A^\alpha) = c_0$.
%By Claim 2, $B^L$ induces in $C$ a structure $B'$ isomorphic to $B$.

Since $|L|=m$ and since the embeddings
$g_\alpha$, $g_\beta$ are distinct whenever
$\alpha,\beta$ are distinct elements of $L$,
we conclude that all of the $m$ 
copies of $\bA$ in the
structure induced by $h(B)$ in $\bC$
have the same colour under $\chi$, which concludes the proof. 
\end{proof}

To finally prove Theorem~\ref{thm:ordered-graph-ramsey},
we combine the partite lemma (Lemma~\ref{lem:partite}) with the so-called 
\emph{partite construction}; again, 
we follow~\cite{NesetrilRoedl}.
%, and again we give more detail
%here since the result will be generalised later. 

\begin{proof}[Proof of Theorem~\ref{thm:ordered-graph-ramsey}]
Let $\bA,\bB$ be $\tau \cup \{\preceq\}$-structures where $\preceq$ denotes a linear order, 
and $r \in \mN$ be arbitrary.
Set $a := |A|$ and $b := |B|$. 
We view $\bA$ as an $a$-partite transversal and $\bB$ as a $b$-partite transversal. 
Let $p \in \mN$ be such that
$([p],<) \to ([b],<)^{([a],<)}_r$ which exists since
$\LO$ is a Ramsey class (Example~\ref{expl:lo}). 
Let $q := {p \choose q}$, 
%$q := |{([p],<) \choose ([a],<)}|$, 
and 
${([p],<) \choose ([a],<)} = \{g_1,\dots,g_q\}$. Construct $p$-partite $\tau \cup \{\preceq\}$-structures $\bP_0,\bP_1,\dots,\bP_q$ inductively as follows. 
Let $\bP_0$ be
such that for any 
$b$ parts $P_{0,i_1},\dots,P_{0,i_b}$
of $\bP_0$ there is an embedding of $\bB$ into
the substructure of $\bP_0$ induced by those parts.
It is clear that such a $(\tau \cup \{\preceq\})$-structure $\bP_0$ exists; one may for instance take an appropriate quasi-ordering $\preceq$ on a disjoint union of the $\tau$-reduct of $\bB$. 
%The reduct of such a structure $\bP_0$ without the order can be obtained as a disjoint union
%of copies of $\bB$ without the order; it is straightforward how to define $\prec$ on the resulting
%structure to obtain a structure $\bP_0$ with the requirements above.  Hmh, gefaellt mir nicht. 

Now suppose that we have already constructed the $p$-partite structure $\bP_{k-1}$,
with parts $P_{k-1,1},\dots,P_{k-1,p}$; 
to construct $\bP_{k}$, let $\bD_{k-1}$ be the $a$-partite system induced in $\bP_{k-1}$ by 
$\bigcup_{i \in [a]} P_{k-1,g_k(i)}$. By the partite lemma (Lemma~\ref{lem:partite}) there exists an $a$-partite 
structure $\bE_k$ such that $\bE_k \to (\bD_{k-1})^\bA_r$. 
We construct the $p$-partite structure $\bP_k$ by amalgamating $\bE_k$ with $\bP_{k-1}$ over $\bD_{k-1}$,
for each occurrence of $\bD_{k-1}$ in $\bE_k$.  

Finally, let $\bC$ be the structure obtained from $\bP_q$
by replacing the linear quasi-order $\preceq$ by a (total) linear extension. 
We claim that $\bC \rightarrow (\bB)^\bA_r$. 
Let $\chi \colon {\bC \choose \bA} \to [r]$ be arbitrary. 
%For $k \in [q]$ we will construct embeddings 
%$h_k \in {\bP_k \choose \bP_q}$ such that
%for all $l \in \{k,\dots,q\}$ there exists a $c_l \in [r]$
%such that for all $f \in {\bD_l \choose \bA}$ we have
%$\chi(h_k \circ f) = c_l$. TYPE CLASH
For $k \in \{0,\dots,q\}$ and $l \in \{k,\dots,q\}$, we will construct embeddings $h_{l,k} \in {\bP_l \choose \bP_{k}}$ such that for all $m \in \{k,\dots,l\}$
%for all $m \in [q]$ 
%there exists a $c_m \in [r]$
%such that for all $f \in {\bD_m \choose \bA}$ we have $\chi(h_{1,q} \circ f) = c_l$.
\begin{itemize}
\item 
$h_{l,m} \circ h_{m,k} = h_{l,k}$, and
\item $|\chi(h_{q,m} \circ {\bD_m \choose \bA})| \leq 1$.
\end{itemize}
Our construction is by induction on $k$, starting with $k = q$. For $k = l = q$ we can choose $h_{q,q}$ to be
the identity. Now suppose that 
$h_{l',k'}$ has already been defined for all $k'$ such that $k \leq k' \leq l' \leq q$. 
We want to define $h_{k,k-1}$.
Since $\bE_k \rightarrow (\bD_{k-1})^\bA_c$, 
there exists an $e_{k-1} \in {\bE_k \choose \bD_{k-1}}$
such that $|\chi(h_{q,k} \circ e_{k-1} \circ {\bD_{k-1} \choose \bA})| \leq 1$.
By construction of $\bP_k$, 
the embedding 
$e_{k-1}$ can be extended to 
an embedding 
$h_{k,k-1} \in {\bP_k \choose \bP_{k-1}}$. 
For $m \in \{k,\dots,l\}$,
we define $h_{m,k-1} := h_{m,k} \circ h_{k,k-1}$,
completing the inductive construction.

For all $m \in [q]$ 
there exists a $c_m \in [r]$ 
such that
for all $f \in {\bD_m \choose \bA}$ we have
$\chi(h_{q,m} \circ f) = c_m$. 
Define $\xi(g_m) := c_m$. 
Since $([p],<) \to ([b],<)^{([a],<)}_r$, 
there exists an $h \in {([p],<) \choose ([a],<)}$
and $c \in [r]$ such that for all $h' \in {([b],<) \choose ([a],<)}$ we have 
$\xi(h \circ h') = c$. 
By construction of $\bP_0$, 
there exists a $g \in {\bP_0 \choose \bB}$ 
such that 
$g(B) \subseteq \bigcup_{i \in [a]} P_{0,h(i)}$. 
To show the claim it suffices to prove that $\chi(g_{k,0} \circ g \circ {\bB \choose \bA}) \leq 1$. 
Let $g' \in {\bB \choose \bA}$ be arbitrary. 
Note that 
$g_{k,0} \circ g \circ g' \in {\bD_k \choose \bA}$ 
for some $k \in [q]$. Hence, 
$\chi(g_{q,0} \circ g \circ g') = \chi(g_{q,k} \circ g_{k,0} \circ g \circ g') = c$, finishing the proof of the claim.  
\end{proof}

% THIS IS THE PLACE that should fix forbidden cliques
% since there are e.g. triangles in the partite power even
% when B is just an edge: example: see on paper.

% The general idea is that the structure from the partite
% lemma only serves as a "template" how to 
% glue together some parts of an inductively constructed structure
% to Ramseyfy this set of parts.  

\subsection{Irreducible homomorphically forbidden structures}
\label{sect:forb-partite}
For every $n \geq 2$, the class of all ordered \emph{$K_n$-free} graphs is Ramsey. In fact, something more general is true; in order
to state the result in full generality, we need the following concept. 

%When $\Gamma$ is a relational $\tau$-structure with domain $D$, then the 
%\emph{Gaifman graph} of $\Gamma$ is the graph
%$(V;E)$ where $E(x,y)$ holds iff $x \neq y$ 
A structure $\bF$ is called \emph{irreducible} (in the terminology of~\cite{NesetrilRoedlPartite})
if for any pair of distinct elements $x,y \in F$ there
 exists an $R \in \tau$ and $z_1,\dots,z_h \in F$
such that $(z_1,\dots,z_h) \in R^\bF$ and $x,y \in \{z_1,\dots,z_h\}$. It is straightforward to verify
that for a set $\mathcal F$ of irreducible structures with finite relational signature $\tau$,
the class
$\Forb(\mathcal F)$ has (strong) amalgamation,
is closed under substructures, isomorphism, and has the joint embedding property, and therefore
is an amalgamation class. 

\begin{theorem}[\Nesetril-R\"odl]
\label{thm:irreducible}
Let $\mathcal F$ be a set of finite irreducible $\tau$-structures. Then $\C := \Forb(\mathcal F) * \LO$ is a Ramsey class.
\end{theorem}

\begin{figure}
  	\includegraphics[scale=0.6]{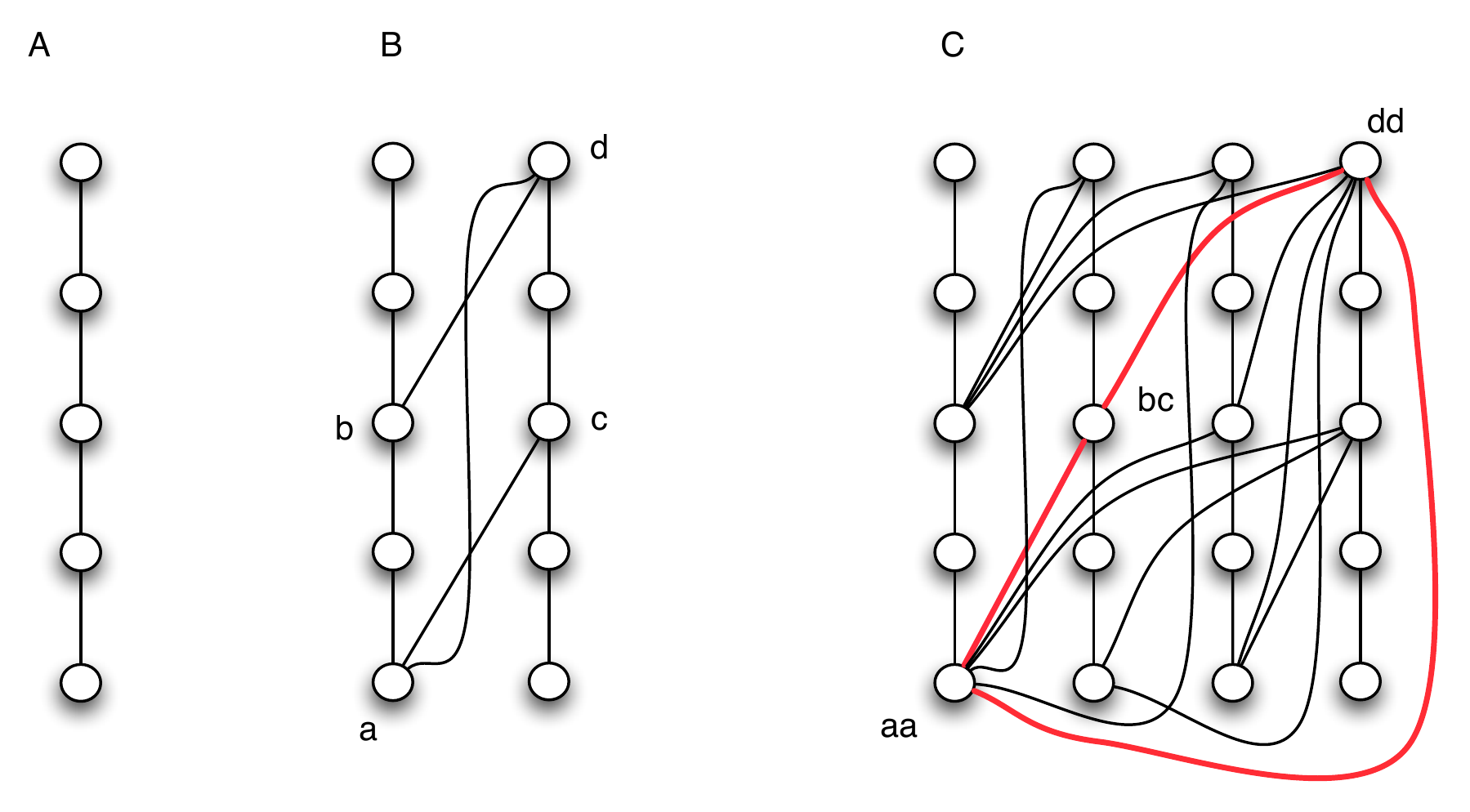} 
\caption{The partite lemma (Lemma~\ref{lem:partite}) can create triangles from triangle-free 5-partite $A$ and $B$.}
\label{fig:partite-counter-ex}
\end{figure}

This theorem can be shown by a variant of the partite method as presented in the previous section. 
However, it is important to note that the proof
from the previous section cannot be applied without an important modification. 
More concretely,
already for the class of triangle-free graphs,
%$n$-partite structure $\bC$ constructed
%in Lemma~\ref{lem:partite} 
NR-powers of $\bB$ over $\bA$
might contain triangles
even if the $n$-partite structures 
$\bA$ and $\bB$ are triangle-free; 
see Figure~\ref{fig:partite-counter-ex}.
To overcome this problem, we need the following definition. Let $\bA,\bB$ be two $n$-partite 
$\tau \cup \{\preceq\}$-structures, and suppose that $\bA$ is transversal. 
Recall that for $u \in \bB$, we write $\pi(u)$ 
for the unique element of $\bA$ that lies
on the same level as $u$. 

\begin{Def}
We
say that $\bA$ is a \emph{template for $\bB$} if
for all $R \in \tau$, 
$(b_1,\dots,b_h) \in R^\bB$ implies that $(\pi(b_1),\dots,\pi(b_h)) \in R^\bA$. 
\end{Def}

We state an important property of the
NR-powers of $\bB$ over $\bA$
when $\bA$ is a template for $\bB$.

\begin{lemma}\label{lem:templates}
Let $\bA$ and $\bB$ be $n$-partite structures
such that $\bA$ is transversal and $\bA$ is a template for $\bB$, and let $r \in \mN$. 
%Let $\bC$ be an NR-power of $\bB$ over $\bA$. 
Then every irreducible structure $\bF$ that
homomorphically maps into 
an NR-power of $\bB$ over $\bA$
%$\bC$ 
also 
homomorphically maps into $\bA$.
\end{lemma}
\begin{proof}
Let $\bC$ be the $d$-th NR-power of $\bB$ 
over $\bA$ for some $d \in \mN$. 
Suppose that $e$ is a homomorphism from
$\bF$ to $\bC$. 
Let $(z_1,\dots,z_h) \in R^\bF$. 
Since $(e(z_1),\dots,e(z_h)) \in \bC$
and by the definition of the NR-power of $\bB$ over
$\bA$, there exists a non-empty set 
$P \subseteq [d]$ and $(w_1,\dots,w_h) \in R^\bB$
such that $(e(z_s))_q = w_s$ for all $q \in P$
and $s \in [h]$. Note that $\pi(w_s) = \pi(z_s)$.
Since $\bA$ is a template for
$\bB$, it follows that
$(\pi(w_1),\dots,\pi(w_h)) \in R^\bA$. 
Hence, $\pi \circ e$ is a homomorphism from
$\bF$ to $\bA$. 
\end{proof}

%The following observation is straightforward to prove, but the core property that explains
%why some results and methods are most naturally phrased in terms of homomorphisms,
%and not embeddings. 
%\begin{lemma}
%Let $A_1,A_2,B$ be $\tau$-structures. 
%Then $A_1 \times A_2$ homomorphically
%embeds $B$ if and only if $A_1$ and $A_2$ homomorphically embed $B$.
%\end{lemma}

We can now modify the partite construction from Section~\ref{sect:all-structs} as follows.
% (a similar idea appeared in~\cite{ProemelVoigtShortProof}). 

\begin{proof}[of Theorem~\ref{thm:irreducible}]
Let $\bA,\bB \in \C$ and $r \in \mN$ be arbitrary. 
By Theorem~\ref{thm:ordered-graph-ramsey},
there exists a $\tau \cup \{\preceq\}$-structure 
$\bC$ where $\preceq$ denotes a linear order (but which need not be from $\C$) such that $\bC \to (\bB)^\bA_r$. 
Let $q := |{\bC \choose \bA}|$,
and ${\bC \choose \bA} = \{g_1,\dots,g_q\}$. 
Let $p := |{\bC \choose \bB}|$,
and ${\bC \choose \bB} = \{f_1,\dots,f_p\}$.
Let $\bC_i$ be the substructure of $\bC$ induced
by $f_i(B)$. 
%$\bigcup_{b \in B} C_{f_i(b)}$. 
We inductively construct a sequence of
$|C|$-partite $\tau$-structures $\bP_0,\bP_1,\dots,\bP_q$.
Let $\bP_0$ be the $(\tau \cup \{\preceq\})$-structure obtained as follows:
define the relation $\preceq$ on the disjoint
union of all the $\bC_i$ by setting $x \preceq y$
if $x$ is a copy of a vertex $x'$
in $\bC$, $y$ is a copy of a vertex $y'$ in $\bC$,
and $x' \preceq y'$ in $\bC$. 
Note that $\bC$ is a template for $\bP_0$. 

The construction of $\bP_k$ for $k > 0$ 
is as in the partite
construction in the proof of Theorem~\ref{thm:ordered-graph-ramsey}:
suppose that we have already constructed
$\bP_{k-1}$; to construct $\bP_k$, let 
$\bD_{k-1}$ be the $|A|$-partite system induced
in $\bP_{k-1}$ by $\bigcup_{i \in [a]} P_{k,g_k(i)}$.
By the partite lemma (Lemma~\ref{lem:partite})
there exists an $|A|$-partite structure $\bE_k$
such that $\bE_k \to (\bD_{k-1})^\bA_r$.
Note that $\bA$ is a template for $\bD_{k-1}$,
and hence, by Lemma~\ref{lem:templates}, 
none of the structures from $\mathcal F$ embeds
into $\bE_k$. We construct the $p$-partite 
structure $\bP_k$ by amalgamating $\bE_k$ with $\bP_{k-1}$ over $\bD_{k-1}$, for each occurrence of $\bD_{k-1}$ in $\bE_k$. 
The proof that $\bP_q \to (\bB)^\bA_r$ is as in 
the proof of Theorem~\ref{thm:ordered-graph-ramsey}. 
\end{proof}

For a recent application of the partite method to prove 
the Ramsey property for classes of structures given 
by homomorphically forbidden trees, see~\cite{Foniok-Trees}.

\section{An inductive proof}
\label{sect:inductive}
In this section we present
a Ramsey class with finite relational signature for
which (to the best of my knowledge) no proof
with the partite method is known. 
Recall the definition of C-relations,
and of convex linear orders
of $C$-relations
from Example~\ref{expl:c-rel}.

\begin{theorem}[see~\cite{BodirskyPiguet,Lee-vorlesungen-ueber-pascaltheorie,Mil79}]
\label{thm:tramsey}
The class 
of all finite binary branching
convexly ordered $C$-relations
is a Ramsey class.
\end{theorem}

This is a consequence of a more powerful theorem due to Milliken~\cite{Mil79}, and follows also from results of Leeb~\cite{Lee-vorlesungen-ueber-pascaltheorie}. A weaker version of this theorem has been shown 
by Deuber~\cite{DeuberTreeRamsey} (my
academic grand-father). 
A direct proof for the statement in the above form
can be found in~\cite{BodirskyPiguet}.

Throughout this section, $\mathcal C$ denotes the class of all finite binary branching convexly ordered $C$-relations. Recall that the members of $\mathcal C$ are in one-to-one correspondence to rooted binary trees, and in the proof it will be convenient to use this perspective. 
%We denote by $\bT(h)$ the $C$-relation of the rooted binary tree with $2^h$ leaves of height $h$, \ie any leaf is at distance $h$ to the root. 
%In particular, $\bT(0)$ denotes the one-element $C$-relation. 

If $\bT$ is a tree with more than one vertex, 
then the root of $\bT$ has
exactly two children; we denote the subtree $\bT$ rooted at the left child (with respect to the convex linear ordering)
by $\bT_{\swarrow}$, and the subtree of $\bT$ rooted at the right child
by $\bT_{\searrow}$ (and we speak of the \emph{left subtree of $\bT$} and the \emph{right subtree of $\bT$}, respectively).
Finally, suppose that $e_1 \in {\bT_\swarrow \choose \bA_{\swarrow}}$
and $e_2 \in {\bT_\searrow \choose \bA_{\searrow}}$,
then $\left< e_1, e_2 \right>$ is the embedding $e$ 
of $\bA$ into $\bT$ defined by $e(a) := e_1(a)$ if $a \in \bA_\swarrow$ and $e(a) := e_2(a)$ if $a \in \bA_\searrow$. 
We write $\bullet$ for the up to isomorphism
unique structure from $\mathcal C$ with one element. 

%For $\bT_1,\bT_2 \in \mathcal C$, define $\bT_1[\bT_2]$
%to be the element of $\mathcal C$ with elements $T_1 \times T_2$ where $C((u_1,v_1);(u_2,v_2),(u_3,v_3))$ holds if
%$C(u_1;u_2,u_3)$ or $u_1 = u_2 = u_3$ and $C(v_1;v_2,v_3)$, and similarly $(u_1,v_1) \prec (u_2,v_2)$
%if $(u_1 \prec u_2)$ or $u_1=u_2$ and $v_1 \prec v_2$. 
%Informally, $\bT_1[\bT_2]$ the the tree obtained from
%$\bT_1$ by replacing every leaf of $\bT_1$ by a copy of $\bT_2$. We also write $\bT^{(1)}$ for $\bT_1$,
%and inductively define $\bT^{(i)} := \bT^{(i-1)}[\bT]$ for $i \geq 2$.
%\begin{lemma}\label{lem:points}
%For all $\bB \in \mathcal C$ and $r \in \mathbb N$ we have
% $\bB^{(r)} \to (\bB)_r^\bullet$. 
%\end{lemma}

%Finally, suppose that $\bT_1$ and $\bT_2$ are disjoint subtrees of $\bT$. Then
%$\left<\bT_1,\bT_2 \right>$ denotes the (uniquely defined) subtree of $\bT$ with leaves
%$T_1 \cup T_2$. Brauchen wir gar nicht. 
%If all copies of $\bA$ in $\bB$ have the same colour, we say that $\bB$ is \emph{$\chi$-monochromatic} (or simply \emph{monochromatic} if the colouring is clear from the context). If the colour is $k$, we also say that $\bB$ is \emph{$k$-chromatic}.

\begin{proof}[of Theorem~\ref{thm:tramsey}]
Let $\bA,\bB \in \mathcal C$, and $r \in \mathbb N$;
we have to show that there is a $\bC \in \mathcal C$
such that $\bC \to (\bB)^\bA_r$.  
We prove the statement by induction over the size of $\bA$.
%, and then over the size of $\bB$. 
For $\bA = \bullet$, 
% the statement follows from Lemma~\ref{lem:points}. 
the proof of the statement 
is easy and left to the reader.  
%We also assume that $\bB$ is distinct from 
%$\bullet$ 
%and that $r \geq 2$ 
%since otherwise there is nothing to be shown. 
%By Lemma~\ref{lem:cover} we can assume
%Similarly as in the proof of Lemma~\ref{lem:partite}, we can assume
%that every element of $\bB$ lies in a copy of $\bA$ in $\bB$.
% TODO: provide the lemma.  
%In particular, $|{\bB \choose \bA}| \geq 1$, 
%Furthermore, $|{\bB_\searrow \choose \bA_\searrow}| \geq 1$
%and $|{\bB_\swarrow \choose \bA_\swarrow}| \geq 1$. 
%The case that $\bB$ has only one element is trivial. 
%Observe that trivially $\bB \rightarrow (\bB)^\bB_r$, and that 
%$\bB \rightarrow (\bB)^\bA_r$ if $|{\bB \choose \bA}|=0$. 
%So assume in the following that $|{\bB \choose \bA}| \geq 1$. 
%Note that $|{\bB_\searrow \choose \bA_\searrow}| \geq 1$
%and $|{\bB_\swarrow \choose \bA_\swarrow}| \geq 1$
%by our assumption that every element of $\bB$ lies
%in a copy of $\bA$ in $\bB$. 

\begin{claim}\label{claim:split}
For all $\bD \in \mathcal C$
there exists an $\bF \in \mathcal C$ such that for any $\chi \colon {\bF \choose \bA} \rightarrow [r]$
 there are $f_1 \in {\bF_\swarrow \choose \bD}$, 
 $f_2 \in {\bF_\searrow \choose \bD}$, and $c \in [r]$ 
such that for all $e_1 \in {\bD \choose \bA_\swarrow}$ and $e_2 \in {\bD \choose \bA_\searrow}$
we have $\chi(\left<f_1 \circ e_1,f_2 \circ e_2\right>) = c$.
\end{claim}
\begin{proof}
By the inductive assumption, there are structures 
$\bF_1,\bF_2 \in \mathcal C$ such that 
$\bF_2 \to (\bD)^{\bA_\searrow}_r$ 
and $\bF_1 \to (\bD)^{\bA_\swarrow}_s$ 
where $s  := \big |[r]^{\bF_2 \choose \bA_\searrow} \big |$. 
Let $\bF$ be such that $\bF_\swarrow = \bF_1$ and
$\bF_\searrow = \bF_2$. 
For a given $\chi \colon {\bF \choose \bA} \to [r]$,
define $\psi \colon {\bF_\swarrow \choose \bA_\swarrow}\rightarrow [r]^{\bF_\searrow \choose \bA_\searrow}$ as follows.
$$ \psi(e_1) := \big ( e_2 \mapsto \chi(\left<e_1,e_2 \right>) \big)$$
By the choice of $\bF_\swarrow = \bF_1$ there exists 
an $f_1 \in {\bF_\swarrow \choose \bD}$ 
and $\phi \colon {\bF_\searrow \choose \bA_\searrow} \to [r]$
such that 
$\psi(f_1 \circ {\bD \choose \bA_\swarrow}) = \{\phi\}$. By the choice of $\bF_\searrow = \bF_2$ there exists an 
$f_2 \in {\bF_2 \choose \bD}$ and a $c \in [r]$
such that $\phi(f_2 \circ {\bD \choose \bA_\searrow}) = \{c\}$. 
Let $g_1 \in {\bD \choose \bA_\swarrow}$
and $g_2 \in {\bD \choose \bA_\searrow}$.
Note that $\psi(f_1 \circ g_1) = \phi$ and
$\phi(f_2 \circ g_2) = c$. 
By definition,  $\phi(f_2 \circ g_2) = 
\chi(\left<f_1 \circ g_1, f_2 \circ g_2 \right >)$,
and hence  $\chi(\left<f_1 \circ g_1, f_2 \circ g_2 \right >) = c$ as desired. 
\end{proof}

% The rest is now like proving the Ramsey property
% for colouring *internal* points in a tree by r colours,
% where we want to find a monochromatic copy of B. 

Let $h$ be the height of $B$ (that is, the maximal distance 
from the root of $\bB$ to one of its leaves), and let $n$ be $h^r$. 
Define $\bC_1,\bC_2,\dots$ inductively as follows. 
Set $\bC_1 := \bullet$, and for $i \geq 2$ let 
$\bC_i$ be the structure $\bF$ that has been constructed for 
$\bD := \bC_{i-1}$ in Claim~\ref{claim:split}. 
Set $\bC := \bC_n$. 

We claim that $\bC \to (\bB)^\bA_r$. 
So let $\chi \colon {\bC \choose \bA} \to [r]$ be given. 
% First the cleaning step. 
For all words $w$ over the alphabet $[2]$ of length $i \in \{0,\dots,n-1\}$ we define $g_w \in {\bC \choose \bC_{n-i}}$, $f_{w1} \in {(\bC_i)_\swarrow \choose \bC_{i-1}}$, $f_{w2} \in {(\bC_i)_\searrow \choose \bC_{i-1}}$, and $c_w \in [r]$ as follows. 
For $i=0$ and $w = \epsilon$, the empty word of length $0$, 
Claim~\ref{claim:split} asserts the existence of $f_1 \in {(\bC_n)_\swarrow \choose \bC_{n-1}}$, $f_2 \in {(\bC_n)_\searrow \choose \bC_{n-1}}$, and $c \in [r]$ such that 
for all 
$e_1 \in {\bC_{n-1} \choose \bA_\swarrow}$ 
and $e_2 \in {\bC_{n-1} \choose \bA_\searrow}$
we have
$\chi(\left<f_1 \circ e_1,f_2 \circ e_2\right>) = c_\epsilon$.
Set $g_1:=f_1$ and $g_2:=f_2$.

Now suppose that $f_w$ and $g_w$ are already defined for a word $w$ of length $i \in \{1,\dots,n-1\}$. 
Let $\psi \colon {\bC_{n-i} \choose \bA} \to [r]$
be the map defined by $\psi(e) := \chi(g_w \circ e)$
for all $e \in {\bC_{n-i} \choose \bA}$.  
Then Claim~\ref{claim:split} asserts the existence of $f_{w1} \in {(\bC_{n-i})_\swarrow \choose \bC_{n-i-1}}$, $f_{w2} \in {(\bC_{n-i})_\searrow \choose \bC_{n-i-1}}$,  and $c_w \in [r]$ such that 
\begin{align}
\psi(\left<f_{w1} \circ e_1,f_{w2} \circ e_2\right>) = c_w
\label{eq:split}
\end{align} 
for all 
$e_1 \in {\bC_{n-i-1} \choose \bA_\swarrow}$ 
and $e_2 \in {\bC_{n-i-1} \choose \bA_\searrow}$. 
Set $g_{w1} := g_w \circ f_{w1}$ 
and $g_{w1} := g_w \circ f_{w2}$.

We claim that there exists an injection
$\beta \colon B \to [2]^{h^r}$ and a $c \in [r]$ 
such that 
\begin{itemize}
\item the map $m$ given by $x \mapsto g_{\beta(x)}(\bullet)$
is from ${\bC \choose \bB}$, and 
\item
for all $b_1,b_2 \in B$ we have that $c_w = c$
when $w$ is the longest common prefix of $\beta(b_1)$
and $\beta(b_2)$. 
\end{itemize}

%Define the colouring $\phi \colon [2]^n \to [r]$ by setting
%$\phi(w) := c_w$. %
%Let $b$ be an element of $\bB$.
%If $b$ has the tree-path $b_1\dots b_{h'}$ for $h' \leq h$ 
%and $b_i \in [2]$,
%let $w_b$ be the 
% word $b_1^{h^{r-1}} \dots b_{h'}^{h^{r-1}} \in %[2]^n$. 

We show this claim by induction on $r$. The statement
is true if $r = 1$ since we can certainly find an
injection $\beta \colon B \to [2]^h$ such that
$x \mapsto g_{\beta(x)}$ is from ${\bC \choose \bB}$, since
$h$ is the height of $\bB$. 

Otherwise we distinguish two possibilities.
We write $S_w$ for the set of words of length at most
$|w|+h^{r-1}$ that start with $w$. Then either 
\begin{enumerate}
\item for every word $w$ of length at most $n' := {n-h^{r-1}} = (h-1) h^{r-1}$ 
there exists a word $u_w \in S_w$ such that $c_{u_w} = r$. In this case, we construct the desired map $\beta$ recursively 
as follows. If $\bB = \bullet$ then define $\beta(b_1) = \epsilon 1\cdots1 \in [2]^n$ (where $\epsilon$ denotes the empty word). 
%All words
%in the image of $\beta$ start with $u := u_\epsilon$.
%If $\bB_\swarrow$ has just one element
%$\{b_1\}$ then define $\beta(b_1) := u1\cdots1 \in [2]^n$.
%If $\bB_\searrow$ has just one element
%$\{b_2\}$ then define $\beta(b_2) := u2\cdots2 \in [2]^n$.
Otherwise, $h \geq 1$, and there exists a word 
$v := u_\epsilon$ with $c_v = r$. 
We repeat this procedure for $\bB_\swarrow$, 
with $v1$ instead of $\epsilon$,
and for $\bB_\searrow$, with $v2$ instead of $\epsilon$,
until $\beta$ is defined on all elements of $\bB$. 
\item there exists a word $w$ of length 
at most $n'$ such that 
$\{c_u \mid u \in S_w\}$ does not contain the colour $r$.
In this case we find an injection 
$\beta \colon B \to S_{w}$ with the desired properties by the inductive hypothesis (since we only consider $r-1$ colours instead of $r$, and $S_{w}$ still has size $h^{r-1}$). 
\end{enumerate}

We claim that $\chi(m \circ e) = c$ 
for all $e \in {\bB \choose \bA}$. 
Since $|A| \geq 2$, 
there are $a_1 \in \bA_\swarrow$
and $a_2 \in \bA_\searrow$. Let $w$ be the longest common prefix of $\beta(e(a_1))$ and 
$\beta(e(a_2))$. 
We then have 
$c_w = c$. 
Write $e$ as $\left<e_1,e_2 \right>$
where $e_1$ is the restriction of $e$ to $\bA_\swarrow$
and $e_2$ is the restriction of $e$ to $\bA_\searrow$. 
%Note that for $i \in [2]$, all words in the image of
%$\beta \circ e_i$ start with $wi$. 
Let $k_1 \in {\bC_{n-i-1} \choose \bA_\swarrow}$ be
$$x \mapsto g_{w1}^{-1}g_{\beta(e_1(x))}(\bullet) \; .$$ 
Similarly, let $k_2 \in {\bC_{n-i-1} \choose \bA_\searrow}$ be $x \mapsto g_{w2}^{-1}g_{\beta(e_2(x))}(\bullet)$. 
% under $
%defined by $x \mapsto g_{\beta(e_1(x))}(\bullet)$,
%and $k_1 \in {\bC \choose \bA_\searrow}$ be
%defined by $x \mapsto g_{\beta(e_2(x))}(\bullet)$.
%Note that $k_i$ can be written as $g_w \circ f_{wi}$ 
We then have $$\chi(m \circ e) = \chi(\left < g_w \circ f_{w1} \circ k_1,g_w \circ f_{w2} \circ k_2 \right >) = c_w = c$$ 
due to Equation~(\ref{eq:split}). 
\end{proof}

% Several problems in the write-up with Diana:
% 1) messed up induction, can't increase the number of
% colours since we only prove the statement for 2 colours. 
% 2) there is an AND statement (see "Why?!" above) where
% only an OR statement is justified
% also typos: 2-chromatic instead of "has just two colours" or so

In Example~\ref{expl:c-rel}, we
have seen that the class of finite ordered binary branching C-relations does not have the Ramsey property. In the context of Conjecture~\ref{conj:ramsey-expansion}, we want to show
how to expand the class to make it Ramsey.

\begin{example}\label{expl:C-with-two-orders}
The class $\mathcal C$ of all finite structures $(L;C,<,\prec)$,
where $<$ is an arbitrary linear order, and $\prec$ is convex with respect to $C$, is a Ramsey class.
This is an immediate consequence of Theorem~\ref{thm:superposition}:
the class $\mathcal C$ can be described as
the superposition of the Ramsey class 
$\LO$ with the class of all convexly ordered $C$-relations, which is Ramsey by Theorem~\ref{thm:tramsey}. 
\end{example}

\section{The ordering property}
\label{sect:ordering}
There are strong links between
the Ramsey property and the
\emph{ordering property} (as defined in~\cite{RamseyClasses,Topo-Dynamics}).

\begin{Def}[Ordering Property]
Let $\mathcal C'$ be a class of finite structures over the signature $\tau \cup \{\prec\}$ where
$\prec$ denotes a linear order, and let $\mathcal C$ be the class of all $\tau$-reducts of structures from $\mathcal C'$. Then $\mathcal C'$ has the \emph{ordering property with respect to $\prec$} if 
for every $\bX \in \mathcal C$
there exists a $\bY \in \mathcal C$ such that
for all expansions $\bX' \in \mathcal C'$ of $\bX$
and $\bY'\in \mathcal C'$ of $\bY$ 
there is an embedding of $\bX'$ into $\bY'$.
\end{Def}

Many examples of classes with the ordering property can
be obtained from Theorem~\ref{thm:ordering} below,
so we rather start with an example of a Ramsey class
\emph{without} the ordering property.

\begin{example}
Let $\mathcal C$ be the class of all finite sets that are
linearly ordered by two linear orders $<_1$ and $<_2$ (see Example~\ref{expl:random-perm}). 
Then $\mathcal C$ does \emph{not} have the ordering property
with respect to $<_1$. Indeed, let $\bA \in \mathcal C$ be the 
structure $(\{0,1,2\};\{(0,1),(1,2),(0,2)\},\{(1,0),(0,2),(1,2)\})$, and
let $\bB$ be an arbitrary $\{<_1\}$-reduct of a structure from $\mathcal C$, that is, an arbitrary finite linearly ordered set. Then 
the expansion of $\bB$ where $<_2$ denotes the same relation as $<_1$ 
is in $\mathcal C$, but certainly contains
no copy of $\bA$. 
%$\{<_2\}$-reduct
%of a structure from $\mathcal C$. 
\end{example}

\begin{prop}\label{prop:ordering}
Let $\Gamma$ be a
homogeneous relational
$\tau$-structure with domain $D$, and suppose that $\Gamma$ has an $\omega$-categorical 
homogeneous expansion $\Gamma'$ with
signature $\tau \cup \{\prec\}$ where $\prec$
denotes a linear order. 
Then the following are equivalent.
\begin{itemize}
\item the age of $\Gamma'$ has the ordering property with respect to $\prec$;
% with a homogeneous
%expansion
%When $\mathcal C$ is the age of a homogeneous
%relational structure $\Gamma = (D;\prec,R_1,R_2,\dots)$, then $\mathcal C$ has the ordering property 
\item for every finite
$X \subseteq D$ there exists a finite $Y \subseteq D$ such that for every $\beta \in \Aut(\Gamma)$ there exists an $\alpha \in \Aut(\Gamma')$ such that $\alpha(X) \subseteq \beta(Y)$. 
\end{itemize}
\end{prop}
\begin{proof}
First suppose that the age $\mathcal C'$ of
$\Gamma'$ 
has the ordering property. Let $X \subset D$ be finite, and let $\bX$ 
be the structure induced by $X$ in $\Gamma$. 
Then there exists a $\bY$ in $\Age(\Gamma)$ such that for all expansions $\bX' \in \mathcal C'$ of $\bX$
and $\bY' \in \mathcal C'$ of $\bY$
there exists an embedding of $\bX'$ into $\bY'$.
%The structure $\bY'$ imbeds into $\mathcal C'$, so 
Suppose without loss of generality that $\bY$ is a substructure of $\Gamma$ with domain $Y$. 
Let $\beta \in \Aut(\Gamma)$ be arbitrary. 
Let $\bX'$ be the structure induced by $\beta(X)$
in $\Gamma'$, and $\bY'$ the structure induced
by $\beta(Y)$ in $\Gamma'$. Since $\beta \in \Aut(\Gamma)$,
$\bX'$ is isomorphic to an expansion of $\bX$,
and $\bY'$ is isomorphic to an expansion of $\bY$.
By assumption, $\bX'$ embeds into $\bY'$. 
 By homogeneity of $\Gamma'$,
this embedding can be extended to an automorphism $\alpha$ of $\Gamma'$,
and $\alpha$ has the desired property. 
 
For the converse,
let $\bX$ be an arbitrary structure in $\Age(\Gamma)$. Let $Z \subseteq D$ be inclusion-wise minimal with the property that for every embedding
$e$ of $\bX$ into $\Gamma$ there
exists an automorphism $\alpha$ of $\Gamma'$
such that $\alpha(e(X)) \subseteq Z$. 
Since $\Gamma'$ is $\omega$-categorical, 
it has a finite number $m$ of orbits of $|X|$-tuples,
and therefore $Z$ has cardinality at most $m |X|$. 
Let $Y \subseteq D$ be such that for every 
$\beta \in \Aut(\Gamma)$ 
there exists an $\alpha \in \Aut(\Gamma')$ 
such that $\alpha(Z) \subseteq \beta(Y)$. 
Let $\bY$ be the structure induced by $Y$ in $\Gamma$. 
Now let $\bX' := (\bX,\prec) \in \mathcal C'$ and 
$\bY' := (\bY,\prec) \in \mathcal C'$ be order expansions of $\bX$ and $\bY$. Let $g$ be
an embedding of $\bY$ into $\Gamma'$. 
%Without loss of generality, $Y'$ is a substructure of $\Gamma'$. 
By the definition of $Z$, there is an
embedding $\rho$ of $\bX'$ into the substructure
induced by $Z$ in $\Gamma'$. 
%By universality of $\Gamma'$,
%there is an isomorphism $\theta$ between
%$Y'$ and a substructure $Y''$ of $\Gamma'$. 
By homogeneity of $\Gamma$, there is
a $\beta \in \Aut(\Gamma)$ that maps $Y$
to $g(Y)$. By the choice of $Y$ 
there exists an $\alpha \in \Aut(\Gamma')$ 
such that $\alpha(Z) \subseteq \beta(Y)$.
Now, $\beta^{-1} \circ \alpha \circ \rho$
is an embedding of $\bX'$ into $\bY'$, which concludes the proof of the ordering property for $\mathcal C'$ with respect to $\prec$. 
\end{proof}

Our next theorem gives a sufficient condition 
for $\omega$-categorical structures to have
the ordering property with respect to a given
ordering; this condition covers
most structures of interest and generalises many 
previous 
isolated results~\cite{Sok-Directed-graphs-and,Proemel-Ramsey-Book,RamseyClasses,Topo-Dynamics}.

An \emph{orbital} of a permutation group $G$ 
on a set $D$ is an orbit $O$ of the componentwise action of $G$ on $D^2$.
An \emph{$O$-cycle} is a sequence of
pairs $(u_1,u_2),(u_2,u_3),\dots,(u_n,u_1)$
from $O$, for some $n$.
We say that $O$ is \emph{cyclic} if 
it contains an $O$-cycle,
and \emph{acyclic} otherwise. 

% Der Fall das \prec bereits fo-defbar in Detla
% ist, wird hiervon abgedeckt: in dem Fall
% splittet natuerlich kein Orbit. 

\begin{theorem}\label{thm:ordering}
Let $\Gamma$ be a 
 homogeneous $\tau$-structure with
domain $D$, and $\prec$ 
an order on $D$
such that $\Gamma' := (\Gamma,\prec)$
is $\omega$-categorical homogeneous Ramsey.
Suppose furthermore that every acyclic orbital of $\Aut(\Gamma)$
is also an orbital of $\Aut(\Gamma')$.
Then $\Age(\Gamma')$ has the ordering property with respect to $\prec$.
\end{theorem}
\begin{proof}
Let $X \subset D$ be finite. By Proposition~\ref{prop:ordering} 
we have to show 
that there exists a finite $Y \subset D$
such that for all $\beta \in \Aut(\Gamma)$ there exists an $\alpha \in \Aut(\Gamma')$ such that 
$\alpha(X) \subseteq \beta(Y)$. 
Since $\Gamma'$ is $\omega$-categorical, 
it has a finite number $m$ of orbits of $|X|$-tuples, and hence there exists a finite $Z \subset D$ with the following properties:
\begin{itemize}
\item for every $\gamma \in \Aut(\Gamma)$ there is a $\delta \in \Aut(\Gamma')$
such that $\delta(\gamma(X)) \subseteq Z$;
\item for every cyclic orbital $O$ of $\Aut(\Gamma)$, $Z$ contains an $O$-cycle. 
\end{itemize}

Since $\Gamma'$ is Ramsey, there
exists by Proposition~\ref{prop:clean}
a finite set $L \subset D$
such that for all 2-element subsets $S_1,\dots,S_\ell$ of $Z$ and all $\chi_i \colon {L \choose S_i} \to [2]$ there exists a $\theta \in \Aut(\Gamma')$
such that $|\chi_i(\theta \circ {Z \choose S_i})| = 1$
for all $i \in [\ell]$.
%such that for all $\eta \in \Aut(\Gamma')$
%and $e_1,e_2 \in {\theta(Z) \choose S}$,
%if $e_1 = \eta \circ e_2$ then
%$\chi(e_1) = \chi(e_2)$. 

Let $\beta \in \Aut(\Gamma)$ be arbitrary.
Define the map $\chi_i \colon {L \choose S_i} \to [2]$ 
as follows. For $g \in {L \choose S_i}$,
put $\chi_i(g) := 0$ if $\beta|_{g(S_i)}$ preserves $\prec$, and $\chi_i(g) := 1$ otherwise. Let $\theta \in \Aut(\Gamma')$ be 
the automorphism that exists for these colourings 
$\chi_1,\dots,\chi_\ell$ according to the choice of $L$. 

We claim that $Y := \theta(Z)$ has the desired
properties, that is, we show that there is 
an $\alpha \in \Aut(\Gamma')$ 
mapping $X$ into $\beta(Y)$. 
By the definition of $Z$, 
there exists a $\delta_1 \in \Aut(\Gamma')$ 
that maps $X$ into $Z$.
%such that $\delta_1(X) \subseteq Z$. 
By the definition of $Z$, 
there also exists a $\delta_2 \in \Aut(\Gamma')$ 
%such that $\delta_2 \beta \theta \delta_1 (X) \subseteq \delta(Z)$. 
that maps $\beta(\theta(\delta_1(X)))$ into $Z$.
%We write $\delta_3$ for $\delta_2 \circ \beta \circ \theta \circ \delta_1$. 

We claim that 
the restriction $g$ of
$\beta \circ \theta \circ \delta_2 \circ \beta \circ \theta \circ \delta_1$
to $X$ can be extended to an automorphism
$\alpha$ of $\Gamma'$. 
Since $\beta,\theta,\delta_1,\delta_2 \in \Aut(\Gamma)$, and by homogeneity of $\Gamma'$,
it suffices to show that $g$ preserves $\prec$. 
So let $x_1,x_2 \in X$ be such that $x_1 \prec x_2$. 
Let $i$ be such that $S_i = \delta_1(\{x_1,x_2\})$, 
and let $T$ be $\theta \circ \delta_1(\{x_1,x_2\})$. 
\begin{itemize}
\item 
If $\chi_i(\theta \circ \delta_1) = 0$,
 then $\beta|_T$ preserves $\prec$. 
It follows that the restriction of $\theta \circ \delta_2 \circ \beta$ to $T$ can be extended to an automorphism $\eta$ of $\Gamma'$. By the property of $\theta$,
this means that 
$\chi_i(\theta \circ \delta_2 \circ \beta \circ \theta \circ \delta_1) = \chi_i(\theta \circ \delta_1) = 0$. By the definition of $\chi_i$, it follows that 
$\beta|_{\theta \circ \delta_2 \circ \beta(T)}$ preserves $\prec$, and so does the restriction of $g$ to $\{x_1,x_2\}$. 
\item 
Otherwise, if $\chi_i(\theta \circ \delta_1) = 1$, then $\beta|_T$ reverses $\prec$. 
In this case the orbital $O$ of $(x_1,x_2)$ in $\Aut(\Gamma)$ cannot be acyclic: it contains 
the orbital $O_1$ of $(x_1,x_2)$ and the orbital $O_2$ of $(\beta(\theta(\delta_1(x_1))),\beta(\theta(\delta_1(x_2))))$ in $\Aut(\Gamma')$, which are distinct, contrary to our assumption for acyclic orbitals. 
Therefore, $Z$ contains an $O$-cycle,
and so does $\theta(Z)$ since $\theta$ preserves $O$. Let 
$(u_0,u_1),(u_1,u_2),\dots,(u_{n-1},u_0)$ be this $O$-cycle in $\theta(Z)$. 
Suppose for contradiction that 
$\chi_i(\theta \circ \delta_2 \circ \beta \circ \theta \circ \delta_1) = 0$. 
We claim that $\beta(u_{i+1}) \prec \beta(u_i)$ 
for all $i \in \{0,\dots,n-1\}$ where the indices are modulo $n$. Then $(\beta(u_{n-1}),\dots,\beta(u_1),\beta(u_0),\beta(u_{n-1}))$ is a directed cycle in $\prec$, a contradiction since $\prec$ is a linear order.
To see the claim, 
observe that if $u_i \prec u_{i+1}$ (this is, $(u_i,u_{i+1}) \in O_1$), then 
$\beta(u_{i+1}) \prec \beta(u_i)$ since $\beta|_T$ reverses $\prec$. On the other hand,
if $u_{i+1} \prec u_i$ (this is, $(u_i,u_{i+1}) \in O_2)$, then $\beta(u_{i+1}) \prec \beta(u_i)$ since $\beta|_{\theta(\delta_2(\beta(T)))}$ preserves $\prec$. 

We conclude that $\chi_i(\theta \circ \delta_2 \circ \beta \circ \theta \circ \delta_1) = 1$,
and thus 
$\beta|_{\theta(\delta_2(\beta(T)))}$ also reverses $\prec$. 
Reversing $\prec$ twice means preserving $\prec$, and 
so we conclude that the restriction of $g = \beta \circ \theta \circ \delta_2 \circ \beta \circ \theta \circ \delta_1$ to $\{x_1,x_2\}$ preserves
$\prec$. 
\end{itemize}
So $g$ indeed preserves $\prec$ on all of $X$,
which proves the claim about the existence of $\alpha \in \Aut(\Gamma')$. 
Note that by the properties of $g$ we also have that $g(X) \subseteq \beta(Y)$, and this concludes the proof. 
\end{proof}

\begin{cor}
The following classes have the ordering property with respect
to $\prec$: 
\begin{itemize}
\item the class of all finite $\prec$-ordered graphs;
\item the class of all $C$-relations over finite
sets which are convexly ordered by $\prec$;
\item the class of all finite $\prec$-ordered directed graphs;
\item the class of all finite partially ordered sets with a linear order $\prec$ that extends the partial order.
\end{itemize}
\end{cor}
\begin{proof}
The \Fresse\ limit $\Gamma$ of the 
first three classes do not have acyclic orbitals.
The Ramsey property for those classes has been 
established earlier in this text, so the statement follows
from Theorem~\ref{thm:ordering}. 

Now, let $\Gamma' = (\Gamma,\prec)$ 
be the \Fresse\ limit of the class
from the last item. A proof of 
the Ramsey property for this class can be found in~\cite{Sokic,SokicIsrael}. 
There is one acyclic orbital in $\Gamma$, 
namely the strict order relation of the poset. 
But since $\prec$ is 
a linear extension of the poset relation, 
this orbital is also an orbital of $\Gamma'$. 
%class of all finite graphs does not have acyclic orbitals. %Similarly,
%the \Fresse\ limit of the class of all $C$-relations 
%over finite sets does not have acyclic orbitals. 
The statement therefore follows again from Theorem~\ref{thm:ordering}. 
\end{proof}

% Converse probably hard. 
% Idea: take acyclic orbit that splits. 
% Take X that contains orbitals from both parts
% of the split. Try to show that for every Y
% we can somewhat order Y so that 
% one of the two parts doesn't appear.
% so X cannot embed into this ordering of Y.

% problem: this approach really has
% something from the reduct classification project
% we take the endo that flips the order some
% where (exists), canonize it modulo constants,
% and try to use it to re-colour Y. 
% this will work in all examples, but in general
% it is hard to do...

The following example shows that the 
sufficient condition for the ordering property
that we gave in Theorem~\ref{thm:ordering}
is not necessary.

\begin{example}
Let $(D;\prec)$ be any countable dense linear order without endpoints. 
By Theorem~\ref{thm:superposition}, the structure
$\Gamma := ({\mathbb Q};\Betw,<) * (D;\prec)$ is Ramsey.
Note that the reduct of $\Gamma$ with signature $\{\Betw,<\}$
is isomorphic to $({\mathbb Q};\Betw,<)$ (this can be shown by a simple back-and-forth argument), so we assume
that $\Gamma$ has domain $\mathbb Q$. 
%, but does
%not have the ordering property. To see this, 
% (unwichtig, oder?)
Observe that 
$<$ is certainly an acyclic orbital in 
$({\mathbb Q};\Betw,<)$,
but it splits into the orbital $\{(x,y) : x \prec y \wedge x < y\}$
and the orbital $\{(x,y) : y \prec x \wedge x < y\}$ of 
$\Aut(\Gamma)$. Hence, the condition 
from Theorem~\ref{thm:ordering} does not apply. 
But nonetheless, the age $\mathcal C$ of 
$({\mathbb Q};\Betw,\prec)$ has the ordering property with respect to $<$. To see this, 
note that for every finite substructure $\bA$ of $({\mathbb Q};\Betw,\prec)$ the only two expansions
of $\bA$ by a linear order such that the expansion is isomorphic
to a structure from $({\mathbb Q};\Betw,\prec,<)$ are $(\bA,<)$ and $(\bA,>)$.   
With this observation it is straightforward to adapt the proof given 
in Theorem~\ref{thm:ordering} to show the ordering property of $\mathcal C$. 
\end{example}

\section{Concluding remarks and open problems}
\label{sect:open}

\subsection{An application}
Ramsey classes are an important tool 
in classifications
of \emph{reducts} of structures. When $\Gamma$ is a structure, a \emph{reduct} of $\Gamma$
is a relational structure $\Delta$ with the same domain as $\Gamma$ such that all the relations of $\Delta$
have a \emph{first-order definition} 
over $\Gamma$, that is, 
for every relation $R$ of $\Delta$ 
there exists a first-order formula $\phi$
over the signature of $\Gamma$  (without parameters)
such that a tuple $t$ is in $R$ if and only $\phi(t)$ holds in $\Gamma$. We say that two reducts
are \emph{interdefinable} if they are reducts of each other. 
Quite surprisingly, countable 
structures $\Gamma$ that are homogeneous in a finite relational language
tend to have finitely many reducts, up to interdefinability (see e.g.~\cite{Bennett-thesis,42,Cameron5,Poset-Reducts,Pon11,RandomReducts,Thomas96}), and Thomas~\cite{RandomReducts} conjectured that
this is always the case. If the age of $\Gamma$ is Ramsey, or a homogeneous expansion of the structure is Ramsey, then this helps in classifying
reducts; we refer to the survey article~\cite{BP-reductsRamsey} for the technical details. 
Note that this application provides 
another motivation to study Conjecture~\ref{conj:ramsey-expansion}: if the conjecture is true,
then this means that Ramsey classes can be used
to attack Thomas' conjecture in general. 

\subsection{Link with topological dynamics}
%The link with topological dynamics is not needed
%in the rest of this article; however, this link 
%has stimulated quite some research activity lately, so we will briefly describe some main lines.
%For detailed definitions, we refer to~\cite{Topo-Dynamics,MuellerPongracz,Zucker14}.  
Whether an amalgamation class $\mathcal C$
has the Ramsey property only depends on
the (topological) automorphism group $\Aut(\Gamma)$ of the \Fresse~limit $\Gamma$ of $\mathcal C$: 
by a theorem of Kechris, Pestov, and Todorcevic~\cite{Topo-Dynamics}, 
the class $\mathcal C$ is Ramsey 
 if and only if $\Aut(\Gamma)$ is
\emph{extremely amenable}, that is, if
every continuous action of $\Aut(\Gamma)$
on a compact Hausdorff space has a fixed point. 
This result has
attracted considerable attention~\cite{GutVan11,MuellerPongracz,Van11,Van12,Zucker13}. 
We would like to mention that recently, 
Melleray, Van Th\'e, and Tsankov~\cite{MelTheTsa} showed that 
a variant of Conjecture~\ref{conj:ramsey-expansion}  
(namely Question~\ref{quest:ramsey-expansion} in Section~\ref{sect:open})
is equivalent
to the
so-called \emph{universal minimal flow}
of $\Aut(\Gamma)$ being metrizable and having
a $G_\delta$ orbit. Even more recently, 
Zucker proved that $\Aut(\Gamma)$
does have a $G_\delta$ orbit~\cite{Zucker14} provided that it is metrisable. Hence,
the task that remains to prove Conjecture~\ref{conj:ramsey-expansion} with the topological
approach is to prove that the universal minimal flow of $\Aut(\Gamma)$ is metrisable (also see Theorem 8.14 in~\cite{Zucker14}). 

These 
developments in topological dynamics are 
promising, but so far every single combinatorial result about
Ramsey classes
that can be proved using topological
dynamics also has a direct combinatorial proof.
The converse is not true: we have seen in this introductory article several combinatorial proofs where no
topological proof is known.

\subsection{Variants of the Ramsey expansion conjecture}
We have seen that many classes of homogeneous structures can be expanded
so that the class of expanded structures becomes Ramsey. Note that if we are
allowed to use any expansion, we can 
trivially turn every class into a Ramsey class,
simply by adding unary predicates such
that for every element of every structure in the
class there is a unary predicate that just contains
this element of the structure. 
This is why it is important to require
in Conjecture~\ref{conj:ramsey-expansion} that the expansion has a finite signature. 

A weaker finiteness condition than being 
homogeneous in a finite relational language
is the requirement that the expansion is 
$\omega$-categorical (recall Theorem~\ref{thm:ryll}). Therefore, a natural variant of Conjecture~\ref{conj:ramsey-expansion} is the following. 

\begin{qn}\label{quest:ramsey-expansion}
Is it true that every $\omega$-categorical structure has an
$\omega$-categorical expansion which is Ramsey? Equivalently, is it true that every closed oligomorphic subgroup of $S_\omega$ has an extremely amenable closed oligomorphic subgroup?
\end{qn}

Formally, Conjecture~\ref{conj:ramsey-expansion} and Question~\ref{quest:ramsey-expansion} are unrelated, since both the hypothesis and the conclusion are stronger in Conjecture~\ref{conj:ramsey-expansion}. A common weakening is the question whether every structure which is homogeneous in a finite relational language has an $\omega$-categorical Ramsey expansion.
Time will show for which set of hypotheses we can obtain which positive results.

%; but to me it appears unlikely that in order to find a Ramsey expansion for a structure which is homogeneous in a finite relational language we might need an infinite signature (typically, an ordering of the domain and first-order definable relations are sufficient). 
% ACHTUNG: doch nicht so klar,
% bedenke dass man mit = Strukturen
% interpretieren kann, die nicht homogen
% in endlicher relationaler Signatur sind!

Also in Question~\ref{quest:ramsey-expansion}, the assumption that
the expansion be $\omega$-categorical is important, 
since otherwise the answer is trivially positive since the trivial group that just consists of the identity element is extremely amenable. 
It is also true (Theorem 4.5 in~\cite{Kechris-nonarchimedian}) that every closed oligomorphic subgroup of $S_\omega$ has a non-trivial extremely amenable subgroup (which is not $\omega$-categorical, though). 

\ignore{
\subsection{Lattices}
Our confidence in a positive answer for Question ~\ref{quest:ramsey-expansion} is smaller than
for Conjecture~\ref{conj:ramsey-expansion}, 
because there are concrete classes of $\omega$-categorical structures (that are \emph{not} homogeneous in a finite relational signature)
where we do not know how to find 
an $\omega$-categorical Ramsey expansion.
% Internal comment: I am also less fluent
% with the proof that vector spaces are Ramsey, 
% it seems that much more than just
% the partite method is needed there. 

An example of such a class is the
class $\mathcal L$ of all finite lattices $(L;\wedge,\vee)$, where $\wedge$ and $\vee$ are binary
function symbols denoting meet and join, respectively. 
It has the amalgamation property (see \cite{GraetzerLattices}), and its \Fresse\ limit is $\omega$-categorical. It is also true that $\mathcal L$ has
the `Ramsey property for points', that is, 
for all $\bB \in \mathcal L$ and $r \in \mN$ there exists
an $\bC \in \mathcal L$ such that 
$\bB \to (\bB)^{\bullet}_r$, where $\bullet$ is
the 1-element lattice~\cite{PosetRamsey}. 
We mention that for the universal homogeneous distributive lattice, an
$\omega$-categorical Ramsey order-expansion
has been found by Miodrag Sokic~\cite{SokicDlattice}. 
}

\subsection{More Ramsey classes from old?}
We do not know the answer to the following question 
about model companions and model-complete
cores in the context of Ramsey classes. 

\begin{qn}\label{quest:mc}
Suppose that $\Gamma$ is a relational structure 
with a homogeneous Ramsey expansion with finite relational signature, and let $\Delta$
be the model companion of $\Gamma$.
Is it true that $\Delta$ has a
homogeneous Ramsey expansion with finite relational signature?
\end{qn}

A positive answer would be a strengthening
of Theorem~\ref{thm:mc-Ramsey}. 
We can ask the same question for
the model-complete core $\Delta$ of $\Gamma$.
Note that a positive answer to Conjecture~\ref{conj:ramsey-expansion} implies a positive 
answer to Question~\ref{quest:mc}. 

We also have a variant for $\omega$-categorical expansions instead of homogeneous expansions in a finite relational language. 

\begin{qn}\label{quest:omega-cat-mc}
Suppose that $\Gamma$ is a relational structure 
with an $\omega$-categorical Ramsey expansion, and let $\Delta$
be the model companion of $\Gamma$.
Is it true that $\Delta$ has an $\omega$-categorical Ramsey expansion?
\end{qn}

%\subsection{On the scope of the partite method}
%Cherlin expressed some reservation to these conjectures since he would prefer a conjecture that would rather suggest how to expand the structure. 
%Weakest version: every homogeneous structure with finite relational signature has
%an $\omega$-categorical Ramsey expansion. 

\ignore{
\subsection{Forbidden thick cycles}
A concrete open problem: describe forbidden thick cycle class and ask for a Ramsey lift. 
% Don't put that: it should be in the 
% scope of the partite method. 
}

\section*{Acknowledgements}
I want to thank Miodrag Sokic for the proof about adding constants,
Diana Piguet for the permission to include parts from our unpublished joint paper, and Lionel Van Th\'e and Miodrag Sokic for discussions about the partite method. I also want to thank
Antoine Mottet and Andr\'as Pongr\'acz for
discussions around the topic of this survey. Many thanks to Trung Van Pham, Andr\'as Pongr\'acz, Miodrag Sokic, and to the anonymous referee for many very helpful comments on earlier versions of the text. 
%(Needless to say that 
%all remaining errors are exclusively my fault.) 

\bibliographystyle{abbrv}
\bibliography{../../global}

\myaddress

\newpage
\thispagestyle{empty}
\mbox{}
\newpage

\end{document}